# Beyond Nonlinear Small-Gain Design: DADS with Partial-State Feedback

**Iasson Karafyllis* and Miroslav Krstic****

*Dept. of Mathematics, National Technical University of Athens,
Zografou Campus, 15780, Athens, Greece,
email: iasonkar@central.ntua.gr

**Dept. of Mechanical and Aerospace Eng., University of California, San Diego, La Jolla, CA 92093-0411, U.S.A., email: krstic@ucsd.edu

*To Eduardo Sontag on the occasion of his 75th birthday*
*with gratitude for his highly influential*
*contribution to control theory*

## Abstract

Eduardo Sontag and coauthors studied Input-to-Output Stability (IOS) and the output asymptotic gain property. These notions changed control theory and recently had an impact on robust adaptive control through the Deadzone-Adapted Disturbance Suppression (DADS) control scheme. Moreover, recently the notion of IOS was extended to systems described by Partial Differential Equations (PDEs). In this work, we celebrate Eduardo Sontag by combining DADS and IOS for PDEs: we study the partial-state regulation problem for a scalar Ordinary Differential Equation (ODE) which is interconnected with a possibly infinite-dimensional system. In such a case the DADS control scheme can allow an escape from the requirements of the small-gain theorem that is mainly used for partial-state feedback. We show the design procedure of partial-state DADS controllers and we prove robust regulation even in the presence of external inputs (disturbances) without assuming knowledge of any disturbance/parameter bounds. The DADS controller is applied to three different cases of the interconnection of an ODE with an almost completely unknown: (a) heat PDE, (b) transport PDE, and (c) wave PDE with viscous damping. We show that the same DADS controller can achieve robust regulation in all three cases.

## 1. Introduction

Robustness in adaptive control is an important issue that has been studied by many researchers in control/stability theory. For time-invariant nonlinear control systems with no persistence of excitation assumed, researchers have used: nonlinear damping (see [20, 21, 10, 12]), leakage (see [2, 3, 27, 30]), supervision for direct adaptive schemes (see [1]), dynamic (high) gains or gain adjustment (see [4, 19, 23, 24, 29]), projection methodologies (see [2] and Appendix E in [20]), and deadzone in the update law (introduced in the paper [28] and well explained in the book [2]). Recently, the Deadzone-Adapted Disturbance Suppression (DADS) adaptive control scheme (see [13, 14, 15, 16, 17]) provided a simple, direct, adaptive control scheme that combines three elements: (a) nonlinear damping (as in [20, 12]), (b) single-gain adjustment (the dynamic feedback has only one state), and (c) deadzone in the update law. The DADS controller achieves attenuation of the plant state to an assignable small level, despite the presence of (time-varying) disturbances and unknown (time-varying) parameters of arbitrary and unknown bounds. Moreover, the DADS controller prevents gain and state drift regardless of the size of the disturbance and unknown parameter. The latter property is a consequence of the practical Uniform Bounded-Input-Bounded-State (p-UBIBS) property and the practical Input-to-Output Stability (p-IOS) property that both hold for the closed-loop system. Therefore, DADS guarantees stability notions that were studied by Eduardo Sontag and coauthors in [34, 35, 36, 37].

Eduardo Sontag changed mathematical control theory by introducing Input-to-State Stability (ISS) and studying Input-to-Output Stability (IOS) in [31, 32, 33, 34, 35, 36, 37]. Recently, ISS and IOS has been extended to the case of systems described by Partial Differential Equations (PDEs): see for instance [8, 26].

This paper is a link between IOS, PDEs and robust adaptive control. We celebrate Eduardo Sontag's 75th birthday by combining DADS and IOS for PDEs: we study the partial-state regulation problem for a scalar Ordinary Differential Equation (ODE) which is interconnected with a possibly infinite-dimensional system. The main tool for robust partial-state feedback design is using in an instrumental way ISS and IOS: it is the small-gain theorem given in [5]. Here we show that the DADS control scheme can allow an escape from the requirements of the small-gain theorem. We show the design procedure of partial-state DADS controllers and we prove robust regulation even in the presence of external inputs (disturbances) without assuming knowledge of any disturbance/parameter bounds. We also show that the DADS feedback law is successful even in the case where the neglected dynamics are infinite-dimensional.

The structure of the paper is as follows. First, we start with a subsection that provides the notation and all the stability notions used in the paper (some of them are modifications of well-known notions presented in [5, 7, 8, 25, 31, 32, 33, 34, 35, 36, 37]). Section 2 studies a simple example and shows how DADS can allow the escape from the requirements of a small-gain-based feedback design. The main result (Theorem 1) is stated and discussed in Section 3. Section 4 is devoted to the application of the DADS controller to three different cases of the interconnection of an ODE with an almost completely unknown: (a) heat PDE, (b) transport PDE, and (c) wave PDE with viscous damping. We show that the same DADS controller can achieve robust regulation in all three cases. Section 5 of the paper contains the proof of the main result. The concluding remarks of the present work are provided in Section 6.

**Notation and Basic Notions.** Throughout this paper, we adopt the following notation.

* $\mathbb{R}_+ := [0,+\infty)$. For a vector $x \in \mathbb{R}^n$, $|x|$ denotes its Euclidean norm and $x'$ denotes its transpose.

* Let $D \subseteq \mathbb{R}^n$ be an open set and let $S \subseteq \mathbb{R}^n$ be a set that satisfies $D \subseteq S \subseteq cl(D)$, where $cl(D)$ is the closure of $D$. By $C^0(S\,;\Omega)$, we denote the class of continuous functions on $S$, which take values in $\Omega \subseteq \mathbb{R}^m$. By $C^k(S\,;\Omega)$, where $k \ge 1$ is an integer, we denote the class of functions on

$S \subseteq \mathbb{R}^n$, which take values in $\Omega \subseteq \mathbb{R}^m$ and have continuous derivatives of order $k$. In other words, the functions of class $C^k(S\,;\Omega)$ are the functions which have continuous derivatives of order $k$ in $D=\mathrm{int}(S)$ that can be continued continuously to all points in $\partial D \cap S$. When $\Omega=\mathbb{R}$ then we write $C^0(S)$ or $C^k(S)$. We also use the notation $C^\infty(S\,;\Omega)=\bigcap_{k\geq 0} C^k(S\,;\Omega)$ for the class of smooth functions on $S \subseteq \mathbb{R}^n$.

* Let $X$ be a given normed linear space with norm $\|\cdot\|_X$, let $B \subseteq X$ be a subset of $X$ and let $I \subseteq \mathbb{R}$ be a non-empty interval. By $C^0(I;B)$ we denote the class of continuous mappings $f: I \to B$, i.e., the class of mappings $f: I \to B$ for which the following property holds: for every $t_0 \in I$ and for every $\varepsilon > 0$ there exists $\delta > 0$ such that $\|f(t)-f(t_0)\|_X < \varepsilon$ for all $t \in I \cap (t_0-\delta, t_0+\delta)$. By $C_c^0(I;X)$ we denote the class of continuous mappings $f: I \to X$ with compact support. By $C_c^1(I)$, where $I \subseteq \mathbb{R}$ is an open interval, we denote the class of continuously differentiable functions $f: I \to \mathbb{R}$ with compact support. By $C^1(I;B)$ we denote the class of continuous mappings $f \in C^0(I;B)$ for which the following property holds: there exists $g \in C^0(I;X)$ such for every $t_0 \in I$ and for every $\varepsilon > 0$ there exists $\delta > 0$ with $\left\| \frac{f(t)-f(t_0)}{t-t_0} - g(t_0) \right\|_X < \varepsilon$ for all $t \in I \cap \big((t_0-\delta, t_0) \cup (t_0, t_0+\delta)\big)$. In this case we write $f_t = g$. A mapping $f: I \to B$ is of class $C^2(I;B)$ if $f \in C^1(I;B)$ and $f_t \in C^1(I;X)$.

* Let $X$ be a given Banach space with norm $\|\cdot\|_X$, $D \subseteq X$ be a subset of $X$ and let $I \subseteq \mathbb{R}$ be a given interval with non-empty interior. A function $f: I \to D$ is measurable if there exists a set $E \subset I$ of measure zero and a sequence $\left\{ f_n \in C_c^0(I;X) : n \geq 1 \right\}$ such that $\lim_{n\to+\infty}\big(f_n(t)\big) = f(t)$ for all $t \in I \setminus E$. For $p \in [1,+\infty)$, $L^p(I;D)$ is the set of equivalence classes of measurable functions $f: I \to D$ with $\|f\|_p := \left( \int_I \|f(t)\|_X^p \, dt \right)^{1/p} < +\infty$. $L^\infty(I;D)$ is the set of equivalence classes of measurable functions $f: I \to D$ with $\|f\|_\infty := \operatorname{ess\,sup}_{t\in I}\big(\|f(t)\|_X\big) < +\infty$. When $D=\mathbb{R}$ then we write $L^p(I)$ and $L^\infty(I)$ instead of $L^p(I;\mathbb{R})$ and $L^\infty(I;\mathbb{R})$.

* Let $a<b$ be given constants and $k \geq 1$ be an integer. For $p \in [1,+\infty]$, $W^{k,p}(a,b)$ denotes the Sobolev space of functions in $L^p(a,b)$ with weak derivatives up to order $k$, all in $L^p(a,b)$. We set $H^k(a,b) = W^{k,2}(a,b)$. The closure of $C_c^1\big((a,b)\big)$ in $H^1(a,b)$ is denoted by $H_0^1(a,b)$.

* Let $I \subseteq \mathbb{R}$ be an interval and let $a<b$ be given constants. Let the function $u: I \times (a,b) \to \mathbb{R}$ be given. We use the notation $u[t]$ to denote the profile of $u$ at certain $t \in I$, i.e., $(u[t])(x) = u(t,x)$ for all $x \in (a,b)$. We also use the notation $u[t]$ when $u: I \to X$ and $X$ is a normed linear space and the notation $u(t)$ when $u: I \to \mathbb{R}^n$, i.e., when $X = \mathbb{R}^n$.

Stability notions for finite-dimensional systems. The stability notions that are described below were studied in [5, 7, 15, 25, 31, 32, 33, 34, 35, 36, 37].

Let $D \subseteq \mathbb{R}^p$ be a given closed set with $0 \in D$, $h: \mathbb{R}^n \to \mathbb{R}^k$ be a continuous map with $h(0) = 0$ and $f: \mathbb{R}^n \times D \to \mathbb{R}^n$ be a locally Lipschitz with respect to $w \in \mathbb{R}^n$ mapping with $f(0,0) = 0$. Consider the control system

$$\dot{w} = f(w,d)\,,\ w \in \mathbb{R}^n\,,\ d \in D \tag{1.1}$$

with output

$$y = h(w) \tag{1.2}$$

We assume that system (1.1) is forward complete, i.e., for every $w_0 \in \mathbb{R}^n$ and for every Lebesgue measurable and locally essentially bounded input $d: \mathbb{R}_+ \to D$ the unique solution $w(t) = \phi(t, w_0; d)$ of the initial-value problem (1.1) with initial condition $w(0) = w_0$ corresponding to input $d: \mathbb{R}_+ \to D$ exists for all $t \geq 0$. We use the notation $y(t, w_0; d) = h(\phi(t, w_0; d))$ for all $t \geq 0$, $w_0 \in \mathbb{R}^n$ and for every Lebesgue measurable and locally essentially bounded input $d: \mathbb{R}_+ \to D$.

We say that system (1.1), (1.2) is *practically Input-to-Output Stable (p-IOS)* if there exist a function $\beta \in KL$, a non-decreasing, continuous function $\gamma: \mathbb{R}_+ \to \mathbb{R}_+$ with $\gamma(0) = 0$ and a constant $\alpha \geq 0$ such that the following estimate holds for all $w_0 \in \mathbb{R}^n$, $t \geq 0$ and for every $d \in L^\infty\left(\mathbb{R}_+; D\right)$:

$$\left|y(t, w_0; d)\right| \leq \beta\left(\left|w_0\right|, t\right) + \gamma\left(\left\|d\right\|_\infty\right) + \alpha \tag{1.3}$$

The constant $\alpha \geq 0$ is called the *residual constant* while the function $\gamma$ is called the *gain function of the input $d \in D$ to the output $y$*. When $\alpha = 0$ then we say that system (1.1), (1.2) is *Input-to-Output Stable (IOS)*.

We say that system (1.1), (1.2) satisfies the *practical Output Asymptotic Gain (p-OAG)* property if there exists a non-decreasing, continuous function $\tilde{\gamma}: \mathbb{R}_+ \to \mathbb{R}_+$ with $\tilde{\gamma}(0) = 0$ and a constant $\tilde{\alpha} \geq 0$ such that the following estimate holds for all $w_0 \in \mathbb{R}^n$ and for every $d \in L^\infty\left(\mathbb{R}_+; D\right)$:

$$\limsup_{t \to +\infty}\left(\left|y(t, w_0; d)\right|\right) \leq \tilde{\gamma}\left(\left\|d\right\|_\infty\right) + \tilde{\alpha} \tag{1.4}$$

The constant $\tilde{\alpha} \geq 0$ is called the *asymptotic residual constant* while the non-decreasing, continuous function $\tilde{\gamma}$ is called the *asymptotic gain function of the input $d \in D$ to the output $y$*. When $\tilde{\gamma} \equiv 0$ we say that system (1.1), (1.2) satisfies the *zero practical Output Asymptotic Gain property (zero p-OAG)*.

When $h(w) = w$ then the word "output" in the above properties is either replaced by the word "state" (e.g., ISS, p-ISS) or is omitted (e.g., p-AG, zero p-AG).

We say that system (1.1) satisfies the *practical Uniform Bounded-Input-Bounded-State (p-UBIBS)* property if there exists a function $\bar{\gamma} \in K_\infty$ and a constant $\bar{\alpha} > 0$ such that the following estimate holds for all $w_0 \in \mathbb{R}^n$ and for every $d \in L^\infty\left(\mathbb{R}_+; D\right)$:

$$\sup_{t\ge 0}\left(\left|\phi(t,w_0;d)\right|\right)\le \bar{\gamma}\left(\left|w_0\right|\right)+\bar{\gamma}\left(\left\|d\right\|_\infty\right)+\bar{\alpha} \tag{1.5}$$

The p-UBIBS property with $\bar{\alpha}=0$ is called the UBIBS property. Clearly, the p-UBIBS property is equivalent to the existence of a continuous function $B:\mathbb{R}^n\times\mathbb{R}_+\to\mathbb{R}_+$ for which the following estimate holds for all $w_0\in\mathbb{R}^n$ and for every $d\in L^\infty\left(\mathbb{R}_+;D\right)$:

$$\sup_{t\ge 0}\left(\left|\phi(t,w_0;d)\right|\right)\le B\left(w_0,\left\|d\right\|_\infty\right) \tag{1.6}$$

Stability notions for infinite-dimensional systems (see [8]): When the state $w$ belongs to a normed linear space $X$ with norm $\|w\|$, the input set $D$ is a subset of a normed linear space $\mathcal{E}$ with norm $\|d\|_{\mathcal{E}}$ and the output map $h:X\to Y$, where $Y$ is a normed linear space with norm $\|y\|_Y$, is a continuous map, then we assume the existence of a set $\Im\subseteq X$ such that for every $w_0\in\Im$ there exists a set of allowable inputs $\wp\left(w_0\right)\subseteq L^\infty\left(\mathbb{R}_+;D\right)$ for which the solution $w[t]$ with initial condition $w[0]=w_0$ that corresponds to arbitrary input $d\in\wp\left(w_0\right)$ has the following properties: (a) the solution exists for all $t\ge 0$ and is unique, (b) for all $\tau\ge 0$ it holds that $w[\tau]\in\Im$ and the shifted input $\left(P_\tau d\right)(s)=d\left(\tau+s\right)$ for $s\ge 0$ is of class $\wp\left(w[\tau]\right)\subseteq L^\infty\left(\mathbb{R}_+;D\right)$. In addition, if we denote the solution $w[t]$ with initial condition $w[0]=w_0$ that corresponds to arbitrary input $d\in\wp\left(w_0\right)$ by $\phi(t,w_0;d)$ then we assume that we also have $\phi\left(t+\tau,w_0;d\right)=\phi\left(t,\phi\left(\tau,w_0;d\right);P_\tau d\right)$ for all $t,\tau\ge 0$, $w_0\in\Im$ and $d\in\wp\left(w_0\right)$. Finally, we assume that $0\in\Im$, $h(0)=0$ and that the zero input $d(t)\equiv 0$ is of class $\wp\left(0\right)$ with $\phi\left(t,0;d\right)=0$ for all $t\ge 0$. In such a case, all the above stability properties are modified in a straightforward way: $\left|y(t,w_0;d)\right|$ and $\left|\phi(t,w_0;d)\right|$ are replaced by $\left\|h\left(\phi(t,w_0;d)\right)\right\|_Y$ and $\left\|\phi(t,w_0;d)\right\|$, respectively.

## 2. Escape from Small-Gain: A Simple Example

In this section we provide a simple example that shows the efficacy of the DADS adaptive scheme. Consider the planar system

$$\begin{aligned}&\dot{w}=-\bar{p}w+\theta_1 y\\&\dot{y}=\theta_2 w+bu+d\\&w,y,u,d\in\mathbb{R}\end{aligned} \tag{2.1}$$

where $\bar{p}>0$ is an unknown constant and $\theta_1,\theta_2\in\mathbb{R}$, $b\in\left(0,+\infty\right)$ are unknown, possibly time-varying parameters. Our aim is to design a partial-state robust adaptive feedback law that depends only on $y\in\mathbb{R}$. In other words, we assume that the component $w\in\mathbb{R}$ of the plant state is not available and it is not measured.

A partial-state, linear, robust feedback law can be designed by using the small-gain theorem (see [5]). Indeed, if $b(t) \geq b_{\min}$, $|\theta_1(t)| \leq \Theta_1$ and $|\theta_2(t)| \leq \Theta_2$ for $t \geq 0$ where $b_{\min}, \Theta_1, \Theta_2 \in (0,+\infty)$ are constants, then the linear feedback law $u = -ky$ achieves the ISS property for the corresponding closed-loop system for

$$k > \frac{\Theta_1 \Theta_2}{\bar{p} b_{\min}} \tag{2.2}$$

This follows by a standard application of the small-gain theorem on system (2.1) with $u = -ky$. Similarly, if we want to achieve

$$\limsup_{t \to +\infty} \left( |y(t)| \right) \leq \sqrt{2\varepsilon} \tag{2.3}$$

for any given $\varepsilon > 0$, then by applying a combination of dynamic gains and leakage, Theorem 4.1 in [4] can give an adaptive feedback law under the following assumptions:

(a) a lower bound for $\bar{p} > 0$ is known,

(b) $\theta_1, \theta_2, d \in \mathbb{R}$ can be time-varying and constants $\Theta_1, \Theta_2, D > 0$ are known such that $|\theta_1(t)| \leq \Theta_1$, $|\theta_2(t)| \leq \Theta_2$ and $|d(t)| \leq D$ for $t \geq 0$,

(c) $b \in (0,+\infty)$ can be time-varying and a constant $b_{\min} > 0$ is known such that $b(t) \geq b_{\min}$ for $t \geq 0$.

More specifically, the statement of Theorem 4.1 in [4] requires that $b \in (0,+\infty)$ is a known constant and $\theta_1 \in \mathbb{R}$ is a constant, but the proof of Theorem 4.1 in [4] can be easily modified to cover the cases described above.

In other words, in order to design a feedback law that achieves (2.3) for any given $\varepsilon > 0$, it is necessary to know the "strength" of the interconnection of the $w-$subsystem $\dot{w} = -\bar{p}w + \theta_1 y$ with the $y-$subsystem $\dot{y} = \theta_2 w + bu + d$. The "strength" of the interconnection is nicely expressed by the ratio $\frac{\Theta_1 \Theta_2}{\bar{p} b_{\min}}$ appearing in the right-hand side of (2.2) and this ratio must be known.

Next, we show how DADS allows the escape from the requirement of knowing the "strength" of the interconnection. Indeed, relations (3.1), (3.2), (3.3), (3.4) (see next section) hold for system (2.1) with

$$X = S = \mathbb{R}, \ \mho = \mathbb{R}, \ p_2 = 1, \ L(w,y) = w, \ \Phi(w) = \frac{1}{\bar{p}} w^2,$$

$$R = K_1 = K_2 = \frac{1}{\bar{p}}, \ \phi(y) \equiv 0, \ G = 0,$$

while $\delta \in \mathcal{E}$ is irrelevant for this example. Indeed, using (2.1) and definition $\Phi(w) = w^2 / \bar{p}$ we get the inequality

$$\dot{\Phi} = -2w^2 + \frac{2}{\bar{p}} \theta_1 w y \leq -w^2 + \frac{1}{\bar{p}^2} \theta_1^2 y^2$$

Let $\varepsilon,\Gamma,\kappa,a,C>0$ be arbitrary constants with $2\kappa>a$. The DADS feedback law is given by the equations

$$u=-\left(\kappa+\exp(z)\right)^{7}\left(K_1+K_2y^2+K_3y^6\right)y \tag{2.4}$$

$$\dot{z}=\Gamma\exp(-z)\left(\frac{1}{2}y^2-\varepsilon\right)^{+} \quad , \quad z\in\mathbb{R} \tag{2.5}$$

where $K_i$, $i=1,2,3$ are constants that satisfy the following inequalities

$$K_1\geq\frac{\kappa+4aC+4\kappa^3+8a\kappa^2}{4a\kappa^6} \tag{2.6}$$

$$K_2\geq\frac{1}{64a^3\kappa^6}\left(\left(\kappa+4aC+4\kappa^3\right)^2+16\kappa^3+64a^2\kappa^4+64a^2\kappa\right) \tag{2.7}$$

$$K_3\geq\frac{\kappa^4+16a^4}{64a^7\kappa^4} \tag{2.8}$$

The controller (2.4), (2.5) is simple: only one integrator is being used. The dynamic gain $z$ increases in order to overcome the effects of $\theta=(\theta_1,\theta_2)\in\mathbb{R}^2$, $b\in(0,+\infty)$ and $d\in\mathbb{R}$. Furthermore, the controller (2.4), (2.5) can guarantee the practical regulation of the plant state and boundedness of solutions for every bounded parameter $\theta=(\theta_1,\theta_2)\in\mathbb{R}^2$, for every bounded disturbances $d\in\mathbb{R}$, $b\in(0,+\infty)$ with $\inf_{t\geq 0}\left(b(t)\right)>0$ and for every initial condition. In order to understand how all these properties are achieved, we provide the following result, which is a direct corollary of Theorem 1 in next section.

**Corollary 1:** *For every $(w_0,y_0,z_0)\in\mathbb{R}^3$ and for every $d\in L^{\infty}(\mathbb{R}_+)$, $\theta_1\in L^{\infty}(\mathbb{R}_+)$, $\theta_2\in L^{\infty}(\mathbb{R}_+)$, $b\in L^{\infty}\left(\mathbb{R}_+;(0,+\infty)\right)$ with $\inf_{t\geq 0}\left(b(t)\right)>0$, the solution of (2.1) with (2.4), (2.5) and initial condition $\left(w(0),y(0),z(0)\right)=(w_0,y_0,z_0)$ satisfies the following estimates for all $t\geq 0$:*

$$2|w(t)|^2+\bar{p}|y(t)|^2\leq\exp(-\mu t)\left(2|w_0|^2+\bar{p}|y_0|^2\right)+\frac{2a\bar{p}\bar{Z}}{\mu\left(\kappa+\exp(z_0)\right)} \tag{2.9}$$

$$z_0\leq z(t)\leq\ln\left(\exp(z_0)+\frac{\Gamma}{4C}|y_0|^2+a\bar{B}\,\frac{2C\left(1+\exp(z_0)\right)+\varepsilon\Gamma}{4C^2\varepsilon\min(1,\kappa)\left(1+\exp(z_0)\right)}\right) \tag{2.10}$$

$$\limsup_{t\to+\infty}\left(|w(t)|\right)\leq\frac{\sqrt{2\varepsilon}}{\bar{p}}\limsup_{t\to+\infty}\left(|\theta_1(t)|\right) \tag{2.11}$$

$$\limsup_{t\to+\infty}\left(|y(t)|\right)\leq\sqrt{2\varepsilon} \tag{2.12}$$

*where*

$$\mu := \min\left( \bar{p}\frac{2\kappa - a}{2\kappa}, 2C \right) \tag{2.13}$$

$$\begin{aligned} \bar{Z} &:= \|d\|_\infty^2 + g\left(\|\theta_2\|_\infty, z_0\right) + g^2\left(\|\theta_2\|_\infty, z_0\right) \\ &+4\inf_{s\geq 0}\left(b(s)\right) g\left( \frac{1}{\inf_{s\geq 0}\left(b(s)\right)}, z_0 \right) + g^2\left( \frac{1}{\bar{p}}\|\theta_1\|_\infty, z_0 \right) \end{aligned} \tag{2.14}$$

$$\begin{aligned} \bar{B} &:= \frac{1}{2}|w_0|^2 + \frac{\bar{p}}{4}|y_0|^2 + \frac{a\bar{p}Z}{2\mu\left(\kappa + \exp\left(z_0\right)\right)} + g\left(\|\theta_2\|_\infty, z_0\right) \\ &+\|d\|_\infty^2 + g^2\left(\|\theta_2\|_\infty, z_0\right) + 2\inf_{s\geq 0}\left(b(s)\right) g\left( \frac{1}{\inf_{s\geq 0}\left(b(s)\right)}, z_0 \right) \end{aligned} \tag{2.15}$$

*and*

$$g\left(s,l\right) := \left( \left(s - \kappa - \exp\left(l\right)\right)^+ \right)^2 \tag{2.16}$$

*Finally, if* $\lim_{t\to+\infty}\left(d(t)\right) = 0$ *and* $\liminf_{t\to+\infty}\left(b(t)\right) \geq 1/\kappa$ *then one of the following holds: either* $\lim_{t\to+\infty}\left(|y(t)|\right) = \lim_{t\to+\infty}\left(|w(t)|\right) = 0$ *or* $\kappa + \exp\left( \lim_{s\to+\infty}\left(z(s)\right) \right) < \limsup_{t\to+\infty}\left( \max\left( \frac{1}{\bar{p}}|\theta_1(t)|, |\theta_2(t)| \right) \right)$.

The success of the DADS feedback law (2.4), (2.5) in tackling this regulation problem without any knowledge of the parameters of the system is evident. The present example shows how we can escape from the small-gain condition with DADS: DADS increases gradually the controller gain to a level where regulation is feasible.

We have applied the controller (2.4) and (2.5) to system (2.1) with constant parameters with values $\theta_1 = 10$, $\theta_2 = 20$, $\bar{p} = 1$ and $b = 0.1$. The parameters for the controllers were

$$\kappa = 2.1,\ \Gamma = 100,\ \varepsilon = 0.00005,\ K_1 = \frac{15}{2}, K_2 = \frac{87}{2}, K_3 = 24$$

for which inequalities (2.6), (2.7), (2.8) hold with $a = 1$ and $C = 80$. The time evolution of the plant states for the solutions of the closed-loop system (2.1), (2.4), (2.5) is shown in Fig. 1 and in Fig.2 for two cases: the disturbance-free case and the persistent disturbance case $d(t) = 3\sin(t)$. The initial conditions are $\left(w_0, y_0, z_0\right) = \left(-0.5, 0.1, -10\right)$. It is clear that in both cases the controller (2.4), (2.5) achieves to bring the plant state to a distance of at most 0.01 from 0.

The time evolution of $z(t)$ and $u(t)$ for the solutions of the closed-loop system (2.1), (2.4), (2.5) is shown in Fig. 3 and Fig.4, respectively. We found that in both cases (the disturbance-free case and the case where $d(t) = 3\sin(t)$) the evolution is similar: $z(t)$ increases and approaches quickly a limit value (-1.326 in the disturbance-free case and -1.26 in the case where $d(t) = 3\sin(t)$). For $t \geq 6$, $z(t)$ remains almost constant with values approximately equal to its limit value. In the persistent disturbance case, the control input oscillates in a way that tries to cancel the effect of the disturbance.

It should also be noticed that Fig. 1 and Fig.2 show that the DADS controller (2.4), (2.5) achieves perfect regulation in the disturbance-free case even in the case where $\liminf_{t \to +\infty} \left( b(t) \right) < 1/\kappa$ and $\kappa + \exp\left( \lim_{s \to +\infty} \left( z(s) \right) \right) < \limsup_{t \to +\infty} \left( \max\left( \frac{1}{\bar{p}} |\theta_1(t)|, |\theta_2(t)| \right) \right)$.

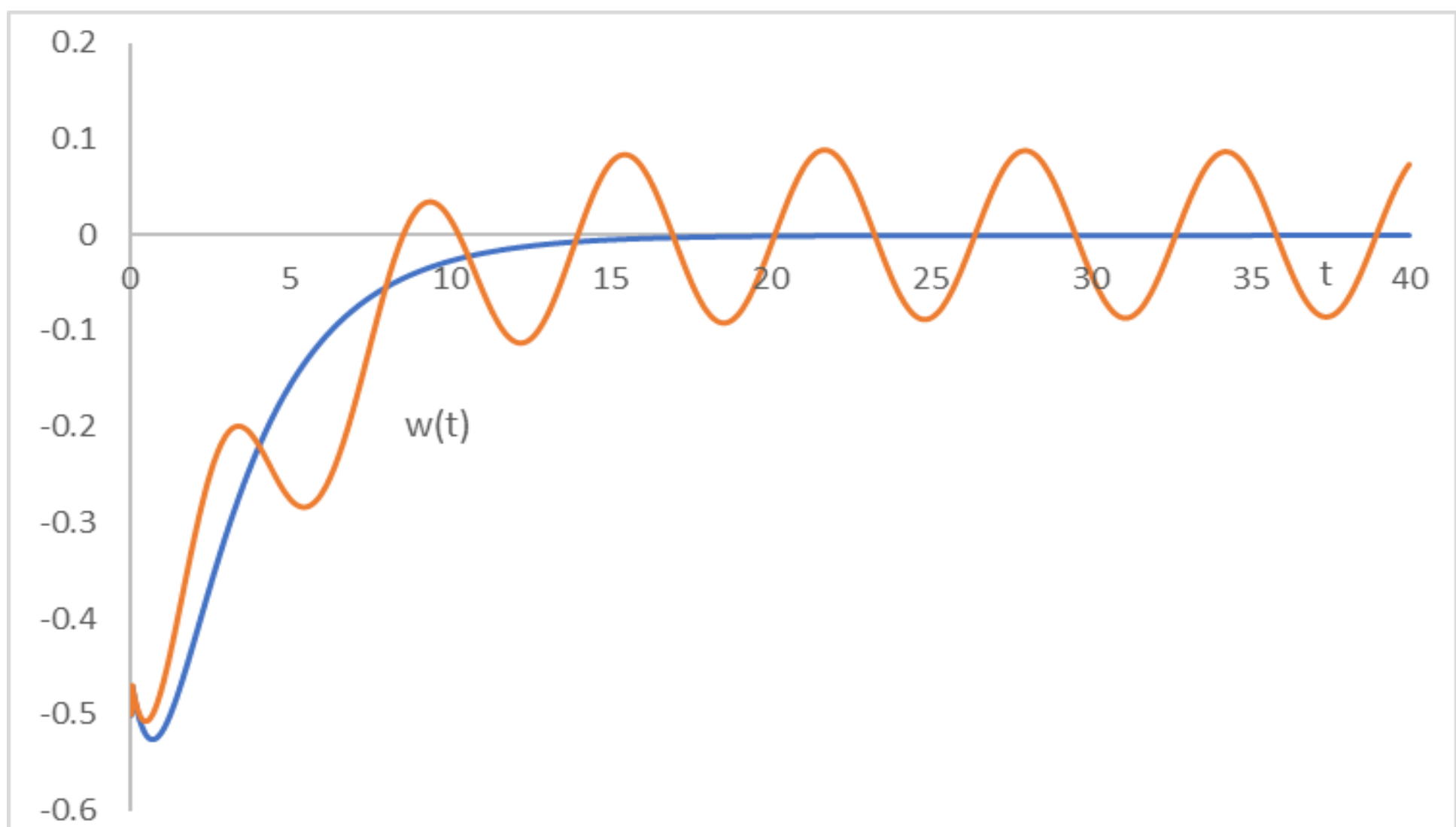

**Fig. 1:** The evolution of the plant state $w$ for the solutions of the closed-loop system (2.1), (2.4), (2.5): blue line the disturbance-free case and red line the persistent disturbance case $d(t) = 3\sin(t)$.

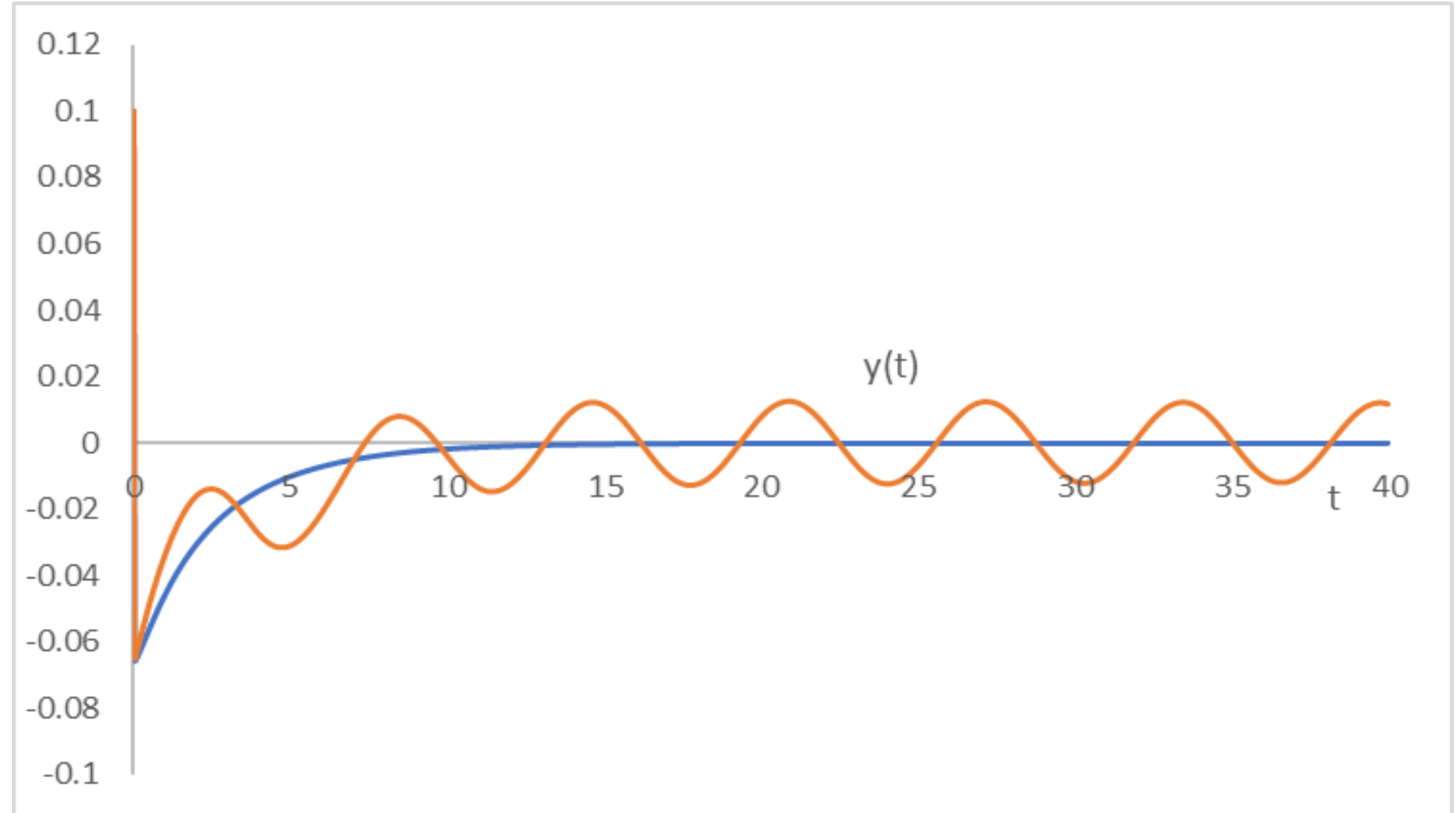

**Fig. 2:** The evolution of the plant state $y$ for the solutions of the closed-loop system (2.1), (2.4), (2.5): blue line the disturbance-free case and red line the persistent disturbance case $d(t) = 3\sin(t)$.

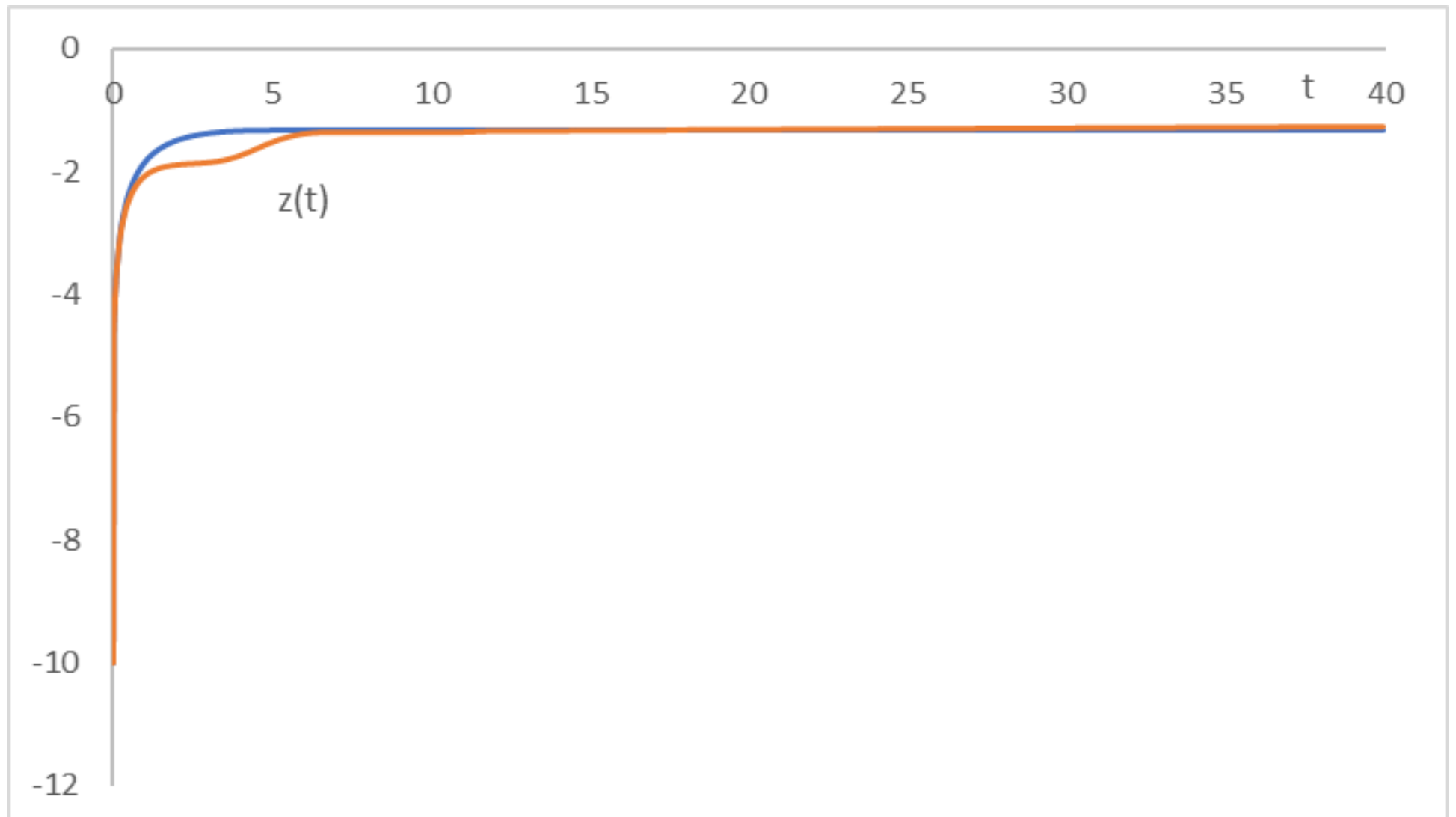

**Fig. 3:** The evolution of the dynamic gain $z$ for the solutions of the closed-loop system (2.1), (2.4), (2.5): blue line the disturbance-free case and red line the persistent disturbance case $d(t) = 3\sin(t)$.

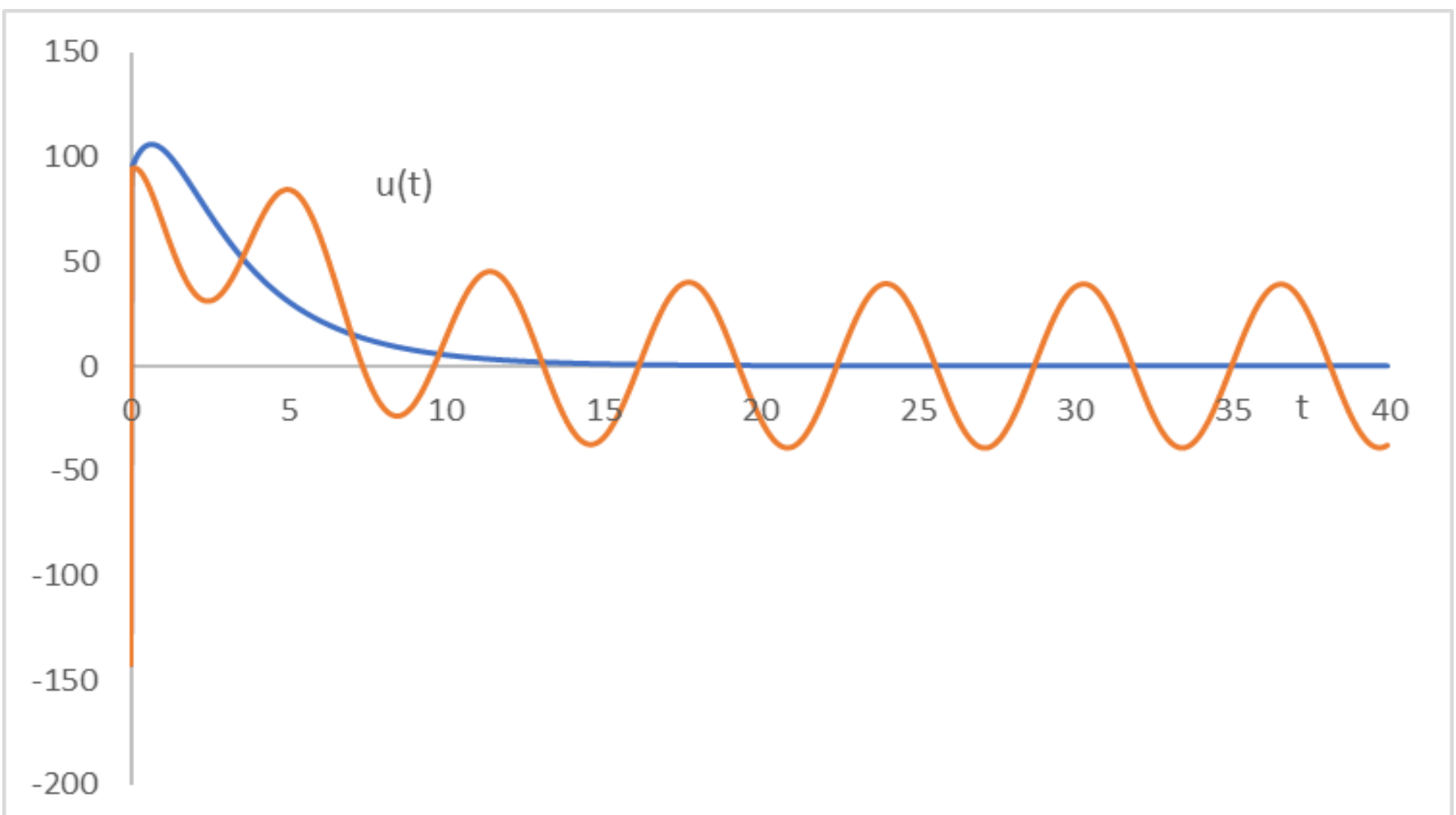

**Fig. 4:** The evolution of the control input $u$ for the solutions of the closed-loop system (2.1), (2.4), (2.5): blue line the disturbance-free case and red line the persistent disturbance case $d(t) = 3\sin(t)$.

## 3. A General Result

In this section we consider the regulation problem of a scalar output $y \in \mathbb{R}$ that is described by the ODE

$$\dot{y} = bu + \left(L(w, y)\right)' \theta_2 + d \tag{3.1}$$

where $u \in \mathbb{R}$ is the control input, $d \in \mathbb{R}$ is a disturbance, $\theta_2 \in \mathbb{R}^{p_2}$, $b > 0$ are unknown possibly time-varying parameters, while $w$ is the unmeasured state component. Here we assume that the state $w$ belongs to a Banach space $X$ with norm $\|w\|$, while $L: S \times \mathbb{R} \to \mathbb{R}^{p_2}$ is an unknown functional where $S \subseteq X$ is a Banach space.

The dynamics of the unmeasured state component $w$ are completely unknown. However, we make the following assumption.

**(A1)** *There exists a continuous functional $\Phi: X \to \mathbb{R}_+$ and constants $K_1, K_2, G, R > 0$ such that*

$$K_1 \|w\|^2 \le \Phi(w) \le K_2 \|w\|^2 \text{, for all } w \in X \tag{3.2}$$

$$\dot{\Phi} \le -\|w\|^2 + R^2 \|\theta_1\|_{\mho}^2 |y|^2 + G \|\delta\|_{\mathcal{E}}^2 \text{, for all } w \in S \tag{3.3}$$

*where $\delta \in \mathcal{E}$, $\theta_1 \in \mho$, where $\mathcal{E}, \mho$ are Banach spaces with norms $\|\theta_1\|_{\mho}$, $\|\delta\|_{\mathcal{E}}$, are unknown possibly time-varying parameters that enter the dynamics of $w$ and $\dot{\Phi}: S \times \mathbb{R} \times \mho \times \mathcal{E} \to \mathbb{R}$ is a functional that is equal to the upper right Dini derivative of $\Phi: X \to \mathbb{R}_+$ for every solution for which $w[t] \in S$, i.e., $\dot{\Phi}(w[t], y(t), \theta_1[t], \delta[t]) = \limsup_{h \to 0^+} \left( h^{-1} \left( \Phi(w[t+h]) - \Phi(w[t]) \right) \right)$ for all $t \ge 0$ for which $w[t] \in S$ and $\limsup_{h \to 0^+} \left( h^{-1} \left( \Phi(w[t+h]) - \Phi(w[t]) \right) \right)$ exists.*

It should be emphasized that none of the constants $K_1, K_2, G, R > 0$ are assumed to be known. Moreover, the functional $\Phi: X \to \mathbb{R}_+$ is not assumed to be known.

We also need an additional assumption for the unknown functional $L: S \times \mathbb{R} \to \mathbb{R}^{p_2}$.

**(A2)** *The functional $L: S \times \mathbb{R} \to \mathbb{R}^{p_2}$ satisfies the following inequality for all $y \in \mathbb{R}$ and $w \in S$*

$$|L(w, y)| \le \|w\| + \phi(y)|y| \tag{3.4}$$

*where $\phi \in C^\infty(\mathbb{R}; \mathbb{R}_+)$ is a known smooth, non-negative function.*

We assume next that $d \in L^\infty(\mathbb{R}_+)$, $\theta_1 \in L^\infty(\mathbb{R}_+; \mho)$, $\theta_2 \in L^\infty(\mathbb{R}_+; \mathbb{R}^{p_2})$, $\delta \in L^\infty(\mathbb{R}_+; \mathcal{E})$, $b \in L^\infty(\mathbb{R}_+; (0, +\infty))$ but we assume no known bounds for $\theta_1 \in \mho$, $\theta_2 \in \mathbb{R}^{p_2}$, $\delta \in \mathcal{E}$, $b \in (0, +\infty)$ and $d \in \mathbb{R}$. We also assume that $\inf_{t \ge 0} (b(t)) > 0$ but we assume no known positive lower bound for $b$.

Our aim is to design a partial-state robust adaptive feedback law that depends only on $y \in \mathbb{R}^n$. In other words, we assume that the component $w \in X$ of the plant state is not available and it is not measured. Here, two important points should be mentioned:

a) the presence of the unknown dynamics that govern the evolution of $w \in X$. That's why in the literature the term "dynamic uncertainty" is used for such systems,

b) the presence of the possibly time-varying unknown input coefficient $b \in (0,+\infty)$ for which nothing is known except the fact that $\inf_{t\geq 0}\left(b(t)\right) > 0$.

Therefore, the control problem studied in the present chapter can equivalently be posed as a partial-state adaptive feedback design problem or as the adaptive feedback design problem for system with dynamic uncertainties.

In order to perform the DADS feedback design, we use the function

$$V(y) = \frac{1}{2}y^2 \tag{3.5}$$

Let $\varepsilon, \Gamma, \kappa, a, C > 0$ be given (arbitrary) constants with $2\kappa > a$. The DADS feedback law is given by the equations

$$u = -(\kappa + \exp(z))^7 \left(P_1(y) + P_2(y)y^2 + P_3(y)y^6\right)y \tag{3.6}$$

$$\dot{z} = \Gamma \exp(-z)\left(V(y) - \varepsilon\right)^+ \quad , \quad z \in \mathbb{R} \tag{3.7}$$

where $P_i \in C^\infty(\mathbb{R})$, $i = 1,2,3$ are smooth functions that satisfy the following inequalities for all $y \in \mathbb{R}$

$$P_1(y) \geq \frac{1}{4a\kappa^6}\left(\kappa + 4aC + 4\kappa^3 + 4a\kappa\phi(y) + 8a\kappa^2\right) \tag{3.8}$$

$$P_2(y) \geq \frac{1}{4a^3\kappa^5}\left(\frac{1}{16\kappa}\left(\kappa + 4aC + 4\kappa^3 + 4a\kappa\phi(y)\right)^2 + \kappa^2 + a^2\phi^2(y) + 4a^2\left(\kappa^3+1\right)\right) \tag{3.9}$$

$$P_3(y) \geq \frac{1}{64a^7\kappa^4}\left(\left(\kappa^2 + a^2\phi^2(y)\right)^2 + 16a^4\right) \tag{3.10}$$

The controller (3.6), (3.7) is a Deadzone-Adapted Disturbance Suppression (DADS) controller, which combines the use of deadzone (in (3.7)) and dynamic nonlinear damping (in (3.6)). The dynamic gain $z$ is being adapted by means of the update law (3.7) and adaptation stops when the output $y$ enters the region defined by $V(y) \leq \varepsilon$. The deadzone prevents the dynamic gain $z$ from growing without bound in the case where a bounded disturbance is present.

The controller (3.6), (3.7) is simple: only one integrator is being used. The dynamic gain $z$ increases in order to overcome the effects of $\theta = (\theta_1, \theta_2) \in \mho \times \mathbb{R}^{p_2}$, $b \in (0,+\infty)$ and $(d,\delta) \in \mathbb{R} \times \mathcal{E}$. Furthermore, the controller (3.6), (3.7) can guarantee the practical regulation of the plant state and boundedness of solutions for every bounded parameter $\theta = (\theta_1, \theta_2) \in \mho \times \mathbb{R}^{p_2}$, for every bounded disturbances $d \in \mathbb{R}$, $\delta \in \mathcal{E}$, $b \in (0,+\infty)$ with $\inf_{t\geq 0}\left(b(t)\right) > 0$ and for every initial condition. In order to understand how all these properties are achieved, we provide the following result.

**Theorem 1:** *Suppose that Assumption (A1) and Assumption (A2) hold. Let $\varepsilon,\Gamma,\kappa,a,C>0$ be constants with $2\kappa>a$. Then every solution*

$$w\in C^0(\mathbb{R}_+;X) \text{ with } \mathbb{R}_+\ni t\to\Phi(w[t]) \text{ being locally absolutely continuous}$$
$$\text{and } w\in L^1((0,T);S) \text{ for each } T>0$$
$$z\in C^1(\mathbb{R}_+),\ y\in C^0(\mathbb{R}_+) \text{ being locally absolutely continuous}$$

*of the closed-loop system (3.1) with (3.6), (3.7) that corresponds to some $d\in L^\infty(\mathbb{R}_+)$, $\theta_1\in L^\infty(\mathbb{R}_+;\mho)$, $\theta_2\in L^\infty(\mathbb{R}_+;\mathbb{R}^{p_2})$, $\delta\in L^\infty(\mathbb{R}_+;\mathcal{E})$, $b\in L^\infty(\mathbb{R}_+;(0,+\infty))$ with $\inf_{t\geq 0}(b(t))>0$ satisfies the following estimates for all $t\geq 0$:*

$$2K_1\|w[t]\|^2+|y(t)|^2\leq \exp(-\mu t)\left(2K_2\|w[0]\|^2+|y(0)|^2\right)+2\mu^{-1}Z \tag{3.11}$$

$$z(0)\leq z(t)\leq \ln\left(\exp(z(0))+\frac{\Gamma}{4C}|y(0)|^2+aB\frac{2C(1+\exp(z(0)))+\varepsilon\Gamma}{4C^2\varepsilon\min(1,\kappa)(1+\exp(z(0)))}\right) \tag{3.12}$$

$$\limsup_{t\to+\infty}(\|w[t]\|)\leq\sqrt{\frac{K_2}{K_1}\left(2\varepsilon R^2\left(\limsup_{t\to+\infty}(\|\theta_1[t]\|_{\mho})\right)^2+G\left(\limsup_{t\to+\infty}(\|\delta[t]\|_{\mathcal{E}})\right)^2\right)} \tag{3.13}$$

$$\limsup_{t\to+\infty}(|y(t)|)\leq\sqrt{2\varepsilon} \tag{3.14}$$

*where*

$$\mu:=\min\left(\frac{2\kappa-a}{2\kappa K_2},2C\right) \tag{3.15}$$

$$\begin{aligned}Z:=\;&G\|\delta\|_\infty^2+\frac{a}{\kappa+\exp(z(0))}\left(\|d\|_\infty^2+g\left(\|\theta_2\|_\infty,z(0)\right)+g^2\left(\|\theta_2\|_\infty,z(0)\right)\right)\\&+\frac{a}{\kappa+\exp(z(0))}\left(4\inf_{s\geq 0}(b(s))\,g\left(\frac{1}{\inf_{s\geq 0}(b(s))},z(0)\right)+g^2\left(R\|\theta_1\|_\infty,z(0)\right)\right)\end{aligned} \tag{3.16}$$

$$\begin{aligned}B:=\;&\frac{1}{2K_1}\left(K_2\|w[0]\|^2+\frac{1}{2}|y(0)|^2+\mu^{-1}Z\right)+g\left(\|\theta_2\|_\infty,z(0)\right)\\&+\|d\|_\infty^2+g^2\left(\|\theta_2\|_\infty,z(0)\right)+2\inf_{s\geq 0}(b(s))\,g\left(\frac{1}{\inf_{s\geq 0}(b(s))},z(0)\right)\end{aligned} \tag{3.17}$$

*and $g$ is given by (2.16). Finally, if $\lim_{t\to+\infty}(d(t))=0$, $\lim_{t\to+\infty}(\delta[t])=0$ and $\liminf_{t\to+\infty}(b(t))\geq 1/\kappa$ then one of the following holds:*

$$\text{either } \lim_{t\to+\infty}(|y(t)|)=\lim_{t\to+\infty}(\|w[t]\|)=0 \text{ or}$$
$$\kappa+\exp\left(\lim_{s\to+\infty}(z(s))\right)<\limsup_{t\to+\infty}\left(\max\left(R\|\theta_1[t]\|_{\mho},|\theta_2(t)|\right)\right) \tag{3.18}$$

It should be noticed that:

- estimate (3.11) shows the p-IOS property of the closed-loop system (3.1) with (3.6), (3.7), output $(w,y)\in X\times\mathbb{R}$ and inputs $\left(\delta,d,\theta_1,\theta_2,\tilde{b}\right)$, where $b=\exp\left(\tilde{b}\right)$,
- the combination of (3.11) and (3.12) shows the p-UBIBS property,
- estimate (3.14) shows the zero p-OAG property for the closed-loop system (3.1) with (3.6), (3.7) and output $y\in\mathbb{R}$. Moreover, estimate (3.14) shows the global practical regulation of the state component $y\in\mathbb{R}$ even when disturbances are present,
- the asymptotic estimates (3.13) and (3.14) show the p-OAG property for the closed-loop system (3.1) with (3.6), (3.7) and output $(w,y)\in X\times\mathbb{R}$,
- exact regulation of the whole state $(w,y)\in X\times\mathbb{R}$ is not excluded. Indeed, if $\kappa+\exp\left(\lim\limits_{s\to+\infty}\left(z(s)\right)\right)\geq\limsup\limits_{t\to+\infty}\left(\max\left(R\left\|\theta_1[t]\right\|_{\mho},\left|\theta_2(t)\right|\right)\right)$ and $\lim\limits_{t\to+\infty}\left(d(t)\right)=0$, $\lim\limits_{t\to+\infty}\left(\delta[t]\right)=0$, $\liminf\limits_{t\to+\infty}\left(b(t)\right)\geq 1/\kappa$, then $\lim\limits_{t\to+\infty}\left(\left|y(t)\right|\right)=\lim\limits_{t\to+\infty}\left(\left\|w[t]\right\|\right)=0$.

Estimates (3.11)-(3.14) and definitions (3.15), (3.16), (3.17), (2.16) indicate how the controller parameters affect the behavior of the closed-loop system. For example, (3.14) shows that the zero p-OAG property holds for the closed-loop system with $y\in\mathbb{R}$ as output and that a smaller $\varepsilon>0$ guarantees a smaller asymptotic residual constant. It also becomes clear from (3.11) -in an almost quantitative way- how the controller parameters affect the gains of the external inputs $d$, $\theta_1,\theta_2,b$, $\delta$. For example, (3.11) shows that an increase of $\kappa$ reduces the gains of $d,b$ and $\theta_1,\theta_2$.

Estimates (3.11), (3.13), (3.14) indicate the three possible phases of convergence that is guaranteed by DADS:

(i) In Phase I, the output $(w,y)$ converges exponentially (with convergence rate $\mu$) in a neighborhood of the set $\left\{(w,y)\in X\times\mathbb{R}:2K_1\left\|w\right\|^2+\left|y\right|^2\leq 2\mu^{-1}Z\right\}$; this is guaranteed by (3.11).

(ii) In Phase II, the output $(w,y)$ converges -in a non-uniform with respect to the initial conditions and inputs way- to a neighborhood of the set $\left\{(w,y)\in X\times\mathbb{R}:2K_1\left|y\right|\leq\sqrt{2\varepsilon},\left\|w\right\|\leq Q\right\}$, where $Q=\sqrt{\frac{K_2}{K_1}\left(2\varepsilon R^2 l_1^2+Gl_\delta^2\right)}$, $l_1=\limsup\limits_{t\to+\infty}\left(\left\|\theta_1[t]\right\|_{\mho}\right)$ and $l_\delta=\limsup\limits_{t\to+\infty}\left(\left\|\delta[t]\right\|_{\mathcal{E}}\right)$; this is guaranteed by (3.13) and (3.14).

(iii) The third phase of convergence may or may nor exist: it is the phase where the output $(w,y)$ converges to zero; this is guaranteed by (3.18) when $\kappa+\exp\left(\lim\limits_{s\to+\infty}\left(z(s)\right)\right)\geq\limsup\limits_{t\to+\infty}\left(\max\left(R\left\|\theta_1[t]\right\|_{\mho},\left|\theta_2(t)\right|\right)\right)$, $\liminf\limits_{t\to+\infty}\left(b(t)\right)\geq 1/\kappa$ and $\lim\limits_{t\to+\infty}\left(d(t)\right)=0$, $\lim\limits_{t\to+\infty}\left(\delta[t]\right)=0$.

## 4. Special Cases

In this section we study three different types of interconnection of a scalar ODE with a PDE. We show that in all three cases the output can be regulated even in the presence of disturbances by the DADS feedback law is given by the equations (3.6), (3.7). Almost nothing is assumed to be known about the PDE itself. Therefore, the DADS feedback law given by the equations (3.6), (3.7) is successful even in the case where the neglected dynamics are infinite-dimensional.

**Case A: Interconnection with a heat PDE.** Consider the following system

$$\begin{aligned} & w_t(t,x) = \bar{p} w_{xx}(t,x) + \theta_1(t,x) K(x, y(t)) + \delta(t,x) \\ & w(t,0) = w(t,1) = 0 \\ & \dot{y}(t) = b(t)u(t) + \theta_2(t) L(w[t], y(t)) + d(t) \end{aligned} \tag{4.1}$$

where $x \in [0,1]$ is the spatial variable, $y(t) \in \mathbb{R}$, $w[t] \in H^2(0,1) \cap H_0^1(0,1)$ are the state components, $u(t) \in \mathbb{R}$ is the control input, $d(t) \in \mathbb{R}$, $\delta[t] \in L^2(0,1)$ are disturbances, $\theta_1[t] \in L^2(0,1)$, $\theta_2(t) \in \mathbb{R}$, $b(t) > 0$ are unknown possibly time-varying parameters, $\bar{p} > 0$ is an unknown constant parameter, $K : [0,1] \times \mathbb{R} \to \mathbb{R}$ is an unknown function and $L : H^2(0,1) \cap H_0^1(0,1) \times \mathbb{R} \to \mathbb{R}$ is an unknown functional that satisfy the following inequalities for all $x \in [0,1]$, $y \in \mathbb{R}$ and $w \in H^2(0,1) \cap H_0^1(0,1)$

$$|L(w,y)| \le \|w\| + \phi(y)|y|, \quad |K(x,y)| \le |y| \tag{4.2}$$

where $\|w\|$ is the norm of $L^2(0,1)$ and $\phi \in C^\infty(\mathbb{R}; \mathbb{R}_+)$ is a known smooth, non-negative function. It should be noted that in order to be able to guarantee existence and uniqueness of solutions of an initial-boundary value problem for system (4.1) additional assumptions are required for the function $K : [0,1] \times \mathbb{R} \to \mathbb{R}$ and the functional $L : H^2(0,1) \cap H_0^1(0,1) \times \mathbb{R} \to \mathbb{R}$. However, next we do not study existence/uniqueness issues for system (4.1) and we focus only on the control issues.

We seek a robust adaptive controller that depends only on $y \in \mathbb{R}$ and guarantees robust practical regulation of the state of system (4.1). The same problem was studied in [16], where the control input coefficient $b(t)$ was assumed to be identically equal to 1, no disturbance $\delta$ appeared in the PDE, $\theta_1$ was spatially constant, i.e., $\theta_1(t,x) = \theta_1(t)$ and $\phi(y)$ was identically equal to zero.

The feedback stabilization problem described above is highly non-trivial. First of all, it should be noticed that almost nothing is known about the dynamics of the unmeasured state component $w$: $\theta_1 \in L^2(0,1)$, $\theta_2 \in \mathbb{R}$, $\bar{p} > 0$ are unknown parameters, $K : [0,1] \times \mathbb{R} \to \mathbb{R}$ and $L : H^2(0,1) \cap H_0^1(0,1) \times \mathbb{R} \to \mathbb{R}$ are unknown functionals that satisfy (4.2). Secondly, it should be emphasized that even when $d = 0$, $\delta = 0$ and $\theta_1(t,x) \equiv \theta_1$, $\theta_2 \in \mathbb{R}$ are constants, the equilibrium point $(y,w) = (0,0)$ of the open-loop system (4.1) can be unstable and system (4.1) can have exponentially increasing (in norm) solutions. For example, when $L(w,y) = -\int_0^1 w(x)dx$, $K(x,y)=\frac{x^2-x-2\bar{p}}{1+2\bar{p}}y$, $6(1+2\bar{p})=\theta_1\theta_2$ and $d=0$, $\delta=0$ the reader can verify that (4.2) is valid and $y(t)=\frac{\theta_2}{6}\exp(t), w(t,x)=\exp(t)x(x-1)$ is an exponentially increasing (in norm) solution of the open-loop system (4.1).

In order to perform the DADS feedback design, we use the functional

$$\Phi(w)=\frac{1}{2\bar{p}}\|w\|^2 \tag{4.3}$$

Using (4.1), (4.2), (4.3), integration by parts, the Cauchy-Schwartz inequality, Wirtinger's inequality ($\|w_x\|^2\geq\pi^2\|w\|^2$ for $w\in H^2(0,1)\cap H_0^1(0,1)$) and the inequalities $\frac{1}{\bar{p}}\|\delta\|\|w\|\leq\frac{4\pi^2-5}{4}\|w\|^2+\frac{1}{\bar{p}^2(4\pi^2-5)}\|\delta\|^2$, $\frac{\|\theta_1\|}{\bar{p}}|y|\|w\|\leq\frac{1}{\bar{p}^2}\|\theta_1\|^2|y|^2+\frac{1}{4}\|w\|^2$ we get for all $w\in H^2(0,1)\cap H_0^1(0,1)$, $y\in\mathbb{R}$, $\theta_1,\delta\in L^2(0,1)$:

$$\begin{aligned}\dot{\Phi}&=-\|w_x\|^2+\frac{1}{\bar{p}}\int_0^1\theta_1(x)K(x,y)w(x)dx+\frac{1}{\bar{p}}\int_0^1\delta(x)w(x)dx\\&\leq-\pi^2\|w\|^2+\frac{\|\theta_1\|}{\bar{p}}|y|\|w\|+\frac{1}{\bar{p}}\|\delta\|\|w\|\\&\leq-\|w\|^2+\frac{1}{\bar{p}^2}\|\theta_1\|^2|y|^2+\frac{1}{\bar{p}^2(4\pi^2-5)}\|\delta\|^2\end{aligned} \tag{4.4}$$

Inequalities (4.2), (4.4) and definition (4.3) show that Theorem 1 can be applied with

$$K_1=K_2=\frac{1}{2\bar{p}},\ X=\mho=\mathcal{E}=L^2(0,1),\ S=H^2(0,1)\cap H_0^1(0,1),$$

$$p_2=1,\ R=\frac{1}{\bar{p}},\ G=\frac{1}{\bar{p}^2(4\pi^2-5)}.$$

Let $\varepsilon,\Gamma,\kappa,a,C>0$ be given (arbitrary) constants with $2\kappa>a$. The DADS feedback law is given by the equations (3.6), (3.7) where $P_i\in C^\infty(\mathbb{R})$, $i=1,2,3$ are smooth functions that satisfy inequalities (3.8), (3.9), (3.10) for all $y\in\mathbb{R}$. Applying Theorem 1 we obtain the following result for system (4.1).

**Corollary 2:** *Suppose that inequality (4.2) holds. Let $\varepsilon,\Gamma,\kappa,a,C>0$ be constants with $2\kappa>a$. Then every solution*

$w\in C^0(\mathbb{R}_+;L^2(0,1))$ *with* $\mathbb{R}_+\ni t\to\|w[t]\|^2$ *being locally absolutely continuous*

*and* $w\in L^1((0,T);H^2(0,1)\cap H_0^1(0,1))$ *for each* $T>0$

$z\in C^1(\mathbb{R}_+), y\in C^0(\mathbb{R}_+)$ *being locally absolutely continuous*

*of the closed-loop system (4.1) with (3.6), (3.7) that corresponds to some* $d \in L^{\infty}\left(\mathbb{R}_{+}\right)$, $\theta_1 \in L^{\infty}\left(\mathbb{R}_{+};L^2\left(0,1\right)\right)$, $\theta_2 \in L^{\infty}\left(\mathbb{R}_{+}\right)$, $\delta \in L^{\infty}\left(\mathbb{R}_{+};L^2\left(0,1\right)\right)$, $b \in L^{\infty}\left(\mathbb{R}_{+};\left(0,+\infty\right)\right)$ *with* $\inf_{t\geq 0}\left(b(t)\right) > 0$, *i.e., for which (4.1), (3.6), (3.7) hold for* $t \geq 0$ *a.e., satisfies estimates (3.11), (3.12), (3.13), (3.14) for all* $t \geq 0$. *Finally, if* $\lim_{t\to+\infty}\left(d(t)\right)=0$, $\lim_{t\to+\infty}\left(\delta[t]\right)=0$ *and* $\liminf_{t\to+\infty}\left(b(t)\right)\geq 1/\kappa$ *then one of the following holds:*

$$\text{either } \lim_{t\to+\infty}\left(\left|y(t)\right|\right) = \lim_{t\to+\infty}\left(\left\|w[t]\right\|\right) = 0 \text{ or}$$

$$\kappa + \exp\left(\lim_{s\to+\infty}\left(z(s)\right)\right) < \limsup_{t\to+\infty}\left(\max\left(\bar{p}^{-1}\left\|\theta_1[t]\right\|, \left|\theta_2(t)\right|\right)\right) \tag{4.5}$$

Similar results can be obtained by applying the same methodology when the boundary conditions in (4.1) are not Dirichlet boundary conditions.

**<u>Case B: Interconnection with a transport PDE.</u>** Consider the following system

$$\begin{aligned} & w_t(t,x) + c w_x(t,x) = \theta_{1,1}\left(t,x\right) K_1(x, y(t)) + \delta\left(t,x\right) \\ & w(t,0) = \theta_{1,2}(t) K_2(y(t)) \\ & \dot{y}(t) = b(t)u(t) + \theta_2(t) L(w[t], y(t)) + d(t) \end{aligned} \tag{4.6}$$

where $x \in \left[0,1\right]$ is the spatial variable, $y(t) \in \mathbb{R}$, $w[t] \in H^1\left(0,1\right)$ are the state components, $u(t) \in \mathbb{R}$ is the control input, $d(t) \in \mathbb{R}$, $\delta[t] \in L^2\left(0,1\right)$ are disturbances, $\theta_1[t] = \left(\theta_{1,1}[t], \theta_{1,2}(t)\right) \in L^2\left(0,1\right) \times \mathbb{R}$, $\theta_2(t) \in \mathbb{R}$, $b(t) > 0$ are unknown possibly time-varying parameters, $c > 0$ is an unknown constant parameter, $K_1 : \left[0,1\right] \times \mathbb{R} \to \mathbb{R}$, $K_2 : \mathbb{R} \to \mathbb{R}$ are unknown functions and $L : H^1\left(0,1\right) \times \mathbb{R} \to \mathbb{R}$ is an unknown functional that satisfy the following inequalities for all $x \in \left[0,1\right]$, $y \in \mathbb{R}$ and $w \in H^1\left(0,1\right)$:

$$\left|L(w,y)\right| \leq \left\|w\right\| + \phi\left(y\right)\left|y\right|, \quad \left|K_1(x,y)\right| \leq \left|y\right|, \quad \left|K_2(y)\right| \leq \left|y\right| \tag{4.7}$$

where $\left\|w\right\|$ is the norm of $L^2\left(0,1\right)$ and $\phi \in C^{\infty}\left(\mathbb{R};\mathbb{R}_{+}\right)$ is a known smooth, non-negative function. Again it should be noted that in order to be able to guarantee existence and uniqueness of solutions of an initial-boundary value problem for system (4.6) additional assumptions are required for the functions $K_1 : \left[0,1\right] \times \mathbb{R} \to \mathbb{R}$, $K_2 : \mathbb{R} \to \mathbb{R}$ and the functional $L : H^1\left(0,1\right) \times \mathbb{R} \to \mathbb{R}$. However, next we do not study existence/uniqueness issues for system (4.6) and we focus only on the control issues.

We seek again a robust adaptive controller that depends only on $y \in \mathbb{R}$ and guarantees robust practical regulation of the state of system (4.6).

It should be emphasized that even when $d = 0$, $\delta = 0$ and $\theta_{1,1}\left(t,x\right) \equiv \theta_{1,1}$, $\theta_{1,2}\left(t\right) \equiv \theta_{1,2}$, $\theta_2\left(t\right) \equiv \theta_2$ are constants, the equilibrium point $\left(y,w\right) = \left(0,0\right)$ of the open-loop system (4.6) can be unstable and system (4.6) can have exponentially increasing (in norm) solutions. For example, when $L(w,y)=\int_0^1 w(x)dx$, $K_1(x,y)=\frac{2(x+c)}{\theta_{1,1}\theta_2}y$, $K_2(y)\equiv 0$, $\theta_{1,1}\theta_2 \geq 2(1+c)$ and $d=0$, $\delta=0$ the reader can verify that (4.7) is valid and $y(t)=\frac{\theta_2}{2}\exp(t), w(t,x)=\exp(t)x$ is an exponentially increasing (in norm) solution of the open-loop system (4.6).

It should also be emphasized that in the case where $K_1(x,y)\equiv 0$, $\delta=0$, $\theta_{1,2}(t)\equiv 1$ and $K_2(y)=y$, the PDE in (4.6) becomes the representation of a state delay, in the sense that its solution satisfies $w(t,x)=y\left(t-c^{-1}x\right)$ for $t>c^{-1}x$. Therefore, regulation problems for systems with state delays like the regulation problem for $\dot{y}(t)=b(t)u(t)+\theta_2(t)\int_{t-c^{-1}}^{t} y(s)ds+d(t)$, where $\theta_2\in L^\infty(\mathbb{R}_+)$, $b\in L^\infty\left(\mathbb{R}_+;(0,+\infty)\right)$ with $\inf_{t\geq 0}(b(t))>0$ and $c>0$ is unknown, can be cast as problems of the form (4.6). It should be noted that such regulation problems cannot be handled in the framework given in [6, 7] without assuming some bounds for $\theta_2$, $\inf_{t\geq 0}(b(t))>0$ and $c>0$.

In order to perform the DADS feedback design, we use the functional the functional:

$$\Phi(w)=\frac{2e}{c}\int_0^1 \exp(-x)w^2(x)dx \tag{4.8}$$

Definition (4.8) implies the following estimates for all $w\in L^2(0,1)$:

$$\frac{2}{c}\|w\|^2 \leq \Phi(w)\leq \frac{2e}{c}\|w\|^2 \tag{4.9}$$

Clearly, definition (4.8), (4.6) and integration by parts give for all $w\in H^1(0,1)$, $y\in\mathbb{R}$, $\delta\in L^2(0,1)$, $\theta_1=(\theta_{1,1},\theta_{1,2})\in L^2(0,1)\times\mathbb{R}$:

$$\begin{aligned}
\dot{\Phi} &= -4e\int_0^1 \exp(-x)w(x)w_x(x)dx+\frac{4e}{c}\int_0^1 \exp(-x)\delta(x)w(x)dx \\
&+\frac{4e}{c}\int_0^1 \exp(-x)\theta_{1,1}(x)w(x)K_1(x,y)dx \\
&=-4e\int_0^1 \exp(-x)\left(\frac{1}{2}w^2(x)\right)_x dx+\frac{4e}{c}\int_0^1 \exp(-x)\delta(x)w(x)dx \\
&+\frac{4e}{c}\int_0^1 \exp(-x)\theta_{1,1}(x)w(x)K_1(x,y)dx \\
&=-2w^2(1)+2ew^2(0)-2e\int_0^1 \exp(-x)w^2(x)dx \\
&+\frac{4e}{c}\int_0^1 \exp(-x)\theta_{1,1}(x)w(x)K_1(x,y)dx+\frac{4e}{c}\int_0^1 \exp(-x)\delta(x)w(x)dx
\end{aligned}$$

Using the Cauchy-Schwartz inequality get for all $w \in H^1(0,1)$, $y \in \mathbb{R}$, $\delta \in L^2(0,1)$, $\theta_1 = (\theta_{1,1}, \theta_{1,2}) \in L^2(0,1) \times \mathbb{R}$:

$$\dot{\Phi} \le -2\|w\|^2 + 2e\theta_{1,2}^2 y^2 + \frac{4e}{c}\|\theta_{1,1}\||y|\|w\| + \frac{4e}{c}\|\delta\|\|w\|$$
$$\le -\|w\|^2 + 2e\left(\theta_{1,2}^2 + \frac{4e}{c}\|\theta_{1,1}\|^2\right)y^2 + \frac{8e^2}{c^2}\|\delta\|^2$$

The above inequality and (4.7), (4.9) show that we can apply Theorem 1 with

$$K_1 = \frac{2}{c},\ K_2 = \frac{2e}{c},\ X = \mathcal{E} = L^2(0,1),\ S = H^1(0,1),$$

$$p_2 = 1,\ R = \sqrt{2e},\ G = \frac{8e^2}{c^2} \text{ and}$$

$$\mho = L^2(0,1) \times \mathbb{R} \text{ with norm } \|\theta_1\|_{\mho} = \|(\theta_{1,1}, \theta_{1,2})\|_{\mho} = \sqrt{\theta_{1,2}^2 + \frac{4e}{c}\|\theta_{1,1}\|^2}$$

Let $\varepsilon, \Gamma, \kappa, a, C > 0$ be given (arbitrary) constants with $2\kappa > a$. The DADS feedback law is given by the equations (3.6), (3.7) where $P_i \in C^\infty(\mathbb{R})$, $i = 1,2,3$ are smooth functions that satisfy inequalities (3.8), (3.9), (3.10) for all $y \in \mathbb{R}$. Applying Theorem 1 we obtain the following result for system (4.6).

**Corollary 3:** *Suppose that inequalities (4.7) hold. Let $\varepsilon, \Gamma, \kappa, a, C > 0$ be constants with $2\kappa > a$. Then every solution*

$$w \in C^0(\mathbb{R}_+; L^2(0,1)) \textit{ with } \mathbb{R}_+ \ni t \to \Phi(w[t]) \textit{ being locally absolutely continuous}$$
$$\textit{and } w \in L^1((0,T); H^1(0,1)) \textit{ for each } T > 0$$
$$z \in C^1(\mathbb{R}_+),\ y \in C^0(\mathbb{R}_+) \textit{ being locally absolutely continuous}$$

*of the closed-loop system (4.6) with (3.6), (3.7) that corresponds to some $d \in L^\infty(\mathbb{R}_+)$, $\theta_1 \in L^\infty(\mathbb{R}_+; L^2(0,1) \times \mathbb{R})$, $\theta_2 \in L^\infty(\mathbb{R}_+)$, $\delta \in L^\infty(\mathbb{R}_+; L^2(0,1))$, $b \in L^\infty(\mathbb{R}_+; (0,+\infty))$ with $\inf_{t \ge 0}(b(t)) > 0$, i.e., for which (4.6), (3.6), (3.7) hold for $t \ge 0$ a.e., satisfies estimates (3.11), (3.12), (3.13), (3.14) for all $t \ge 0$. Finally, if $\lim_{t \to +\infty}(d(t)) = 0$, $\lim_{t \to +\infty}(\delta[t]) = 0$ and $\liminf_{t \to +\infty}(b(t)) \ge 1/\kappa$ then one of the following holds:*

$$\text{either } \lim_{t \to +\infty}(|y(t)|) = \lim_{t \to +\infty}(\|w[t]\|) = 0 \text{ or}$$
$$\kappa + \exp\left(\lim_{s \to +\infty}(z(s))\right) < \limsup_{t \to +\infty}\left(\max\left(\sqrt{2e\theta_{1,2}^2(t) + \frac{8e^2}{c}\|\theta_{1,1}[t]\|^2}, |\theta_2(t)|\right)\right)$$

**Case C: Interconnection with a wave PDE with viscous friction.** Consider the following system

$$\begin{aligned}
&v_{tt}(t,x) = c^2 v_{xx}(t,x) - \sigma v_t(t,x) + \theta_1(t,x) K(x, y(t)) + \delta(t,x) \\
&v(t,0) = v(t,1) = 0 \\
&\dot{y}(t) = b(t)u(t) + \theta_2(t) L(v[t], v_t[t], y(t)) + d(t)
\end{aligned} \tag{4.10}$$

where $x \in [0,1]$ is the spatial variable, $y(t) \in \mathbb{R}$, $v[t], v_t[t]$ are the state components, $u(t) \in \mathbb{R}$ is the control input, $d(t) \in \mathbb{R}$, $\delta[t] \in L^2(0,1)$ are disturbances, $\theta_1[t] \in L^2(0,1)$, $\theta_2(t) \in \mathbb{R}$, $b(t) > 0$ are unknown possibly time-varying parameters, $c, \sigma > 0$ are unknown constant parameters, $K : [0,1] \times \mathbb{R} \to \mathbb{R}$ is an unknown function and $L : \left(H^2(0,1) \cap H_0^1(0,1)\right) \times H_0^1(0,1) \times \mathbb{R} \to \mathbb{R}$ is an unknown functional that satisfy the following inequalities for all $x \in [0,1]$, $y \in \mathbb{R}$, $v \in H^2(0,1) \cap H_0^1(0,1)$ and $v_t \in H_0^1(0,1)$:

$$\left|L(v, v_t, y)\right| \le \sqrt{\|v_x\|_2^2 + \|v_t\|_2^2} + \phi(y)|y|, \quad |K(x,y)| \le |y|, \tag{4.11}$$

where $\|\cdot\|_2$ is the norm of $L^2(0,1)$ and $\phi \in C^\infty(\mathbb{R}; \mathbb{R}_+)$ is a known smooth, non-negative function.

In order to be able to use the framework that is described in the previous section we have to describe the infinite-dimensional part of (4.10) in a different way. Setting $v_t = \varphi$, we can write (4.10) as a first-order in time system:

$$\begin{aligned}
&v_t(t,x) = \varphi(t,x) \\
&\varphi_t(t,x) = c^2 v_{xx}(t,x) - \sigma \varphi(t,x) + \theta_1(t,x) K(x, y(t)) + \delta(t,x) \\
&v(t,0) = v(t,1) = 0 \\
&\dot{y}(t) = b(t)u(t) + \theta_2(t) L(w[t], y(t)) + d(t) \\
&w[t] = \left(v[t], \varphi[t]\right)
\end{aligned} \tag{4.12}$$

The state $w = (v, \varphi)$ is in the Banach space $X = H_0^1(0,1) \times L^2(0,1)$ with norm $\|w\| = \|(v,\varphi)\| = \sqrt{\|v_x\|_2^2 + \|\varphi\|_2^2}$. Again it should be noted that in order to be able to guarantee existence and uniqueness of solutions of an initial-boundary value problem for system (4.12) additional assumptions are required for the functions $K : [0,1] \times \mathbb{R} \to \mathbb{R}$ and the functional $L : \left(H^2(0,1) \cap H_0^1(0,1)\right) \times H_0^1(0,1) \times \mathbb{R} \to \mathbb{R}$. However, next we do not study existence/uniqueness issues for system (4.12) and we focus only on the control issues.

We seek again a robust adaptive controller that depends only on $y \in \mathbb{R}$ and guarantees robust practical regulation of the state of system (4.12).

The constant $\sigma > 0$ is the viscous damping coefficient that provides energy dissipation for the wave PDE in (4.10). Damping mechanisms in wave PDEs were studied in [9, 11, 22]. However, viscous damping does not guarantee bounded solutions for the open-loop system (4.12). Even when $d = 0$, $\delta = 0$ and $\theta_1(t) \equiv \theta_1$, $\theta_2(t) \equiv \theta_2$ are constants, the equilibrium point $(y, w) = (0, 0)$ of the open-loop system (4.12) can be unstable and system (4.12) can have exponentially increasing (in norm) solutions. For example, when $c=1$, $L(v,\varphi,y)=\int_0^1 \sin(\pi x)v(x)dx$, $K_1(x,y)=\sin(\pi x)y$, $\theta_1\theta_2=2\left(1+\sigma+\pi^2\right)$ and $d=0$, $\delta=0$ the reader can verify that (4.11) is valid and $y(t)=\frac{1+\sigma+\pi^2}{\theta_1}\exp(t), v(t,x)=\varphi(t,x)=\exp(t)\sin(\pi x)$ is an exponentially increasing (in norm) solution of the open-loop system (4.12).

We next use the functional defined for all $w=(v,\varphi)\in X$ :

$$\Phi(w)=\frac{(c^2+1)}{\sigma}\|v_x\|_2^2+\frac{1}{\sigma}\|\varphi\|_2^2+\frac{1}{\sigma c^2}\|\varphi+\sigma v\|_2^2 \tag{4.13}$$

Using definition (4.13), Wirtinger's inequality $\pi^2\|v\|_2^2\le\|v_x\|_2^2$, the fact that $\|\varphi+\sigma v\|_2^2\le 2\|\varphi\|_2^2+2\sigma^2\|v\|_2^2$ and the fact that $\|w\|=\|(v,\varphi)\|=\sqrt{\|v_x\|_2^2+\|\varphi\|_2^2}$, we get for all $w=(v,\varphi)\in X$ :

$$\begin{aligned}&\frac{1}{\sigma}\|w\|^2\le\Phi(w)\le\frac{(c^2+1)c^2\pi^2+2\sigma^2}{\sigma c^2\pi^2}\|v_x\|_2^2+\frac{c^2+2}{\sigma c^2}\|\varphi\|_2^2\\&\le\frac{1}{\sigma}\max\left(c^2+1+\frac{2\sigma^2}{c^2\pi^2},1+\frac{2}{c^2}\right)\|w\|^2\end{aligned} \tag{4.14}$$

Using definition (4.13) and (4.12) we get for all $w=(v,\varphi)\in S=\left(H^2(0,1)\cap H_0^1(0,1)\right)\times H_0^1(0,1)$, $y\in\mathbb{R}$, $\delta,\theta_1\in L^2(0,1)$:

$$\begin{aligned}\dot{\Phi}(w)=&\frac{2(c^2+1)}{\sigma}\int_0^1\left(v_x(x)\varphi(x)\right)_x dx-2\|\varphi\|_2^2+2\int_0^1 v(x)v_{xx}(x)dx\\&+\frac{2}{\sigma c^2}\int_0^1\left((c^2+1)\varphi(x)+\sigma v(x)\right)\left(\theta_1(x)K(x,y)+\delta(x)\right)dx\end{aligned} \tag{4.15}$$

Since $\varphi(0)=\varphi(1)=0$, we get from (4.15) and the Cauchy-Schwartz inequality for all $w=(v,\varphi)\in S=\left(H^2(0,1)\cap H_0^1(0,1)\right)\times H_0^1(0,1)$, $y\in\mathbb{R}$, $\delta,\theta_1\in L^2(0,1)$:

$$\begin{aligned}\dot{\Phi}(w)\le&-2\|\varphi\|_2^2-2\|v_x\|_2^2+\frac{2(c^2+1)}{\sigma c^2}\|\varphi\|_2\|\theta_1\|_2|y|\\&+\frac{2(c^2+1)}{\sigma c^2}\|\varphi\|_2\|\delta\|_2+\frac{2}{c^2}\|v\|_2\|\theta_1\|_2|y|+\frac{2}{c^2}\|v\|_2\|\delta\|_2\end{aligned} \tag{4.16}$$

Using the inequalities

$$\|\varphi\|_2\|\theta_1\|_2|y|\le\frac{c^2\sigma}{4(c^2+1)}\|\varphi\|_2^2+\frac{c^2+1}{c^2\sigma}\|\theta_1\|_2^2|y|^2,$$

$$\|\varphi\|_2\|\delta\|_2\le\frac{c^2\sigma}{4(c^2+1)}\|\varphi\|_2^2+\frac{c^2+1}{c^2\sigma}\|\delta\|_2^2,$$

$$\|v\|_2\,\|\theta_1\|_2\,|y| \le \frac{c^2\pi^2}{4}\|v\|_2^2 + \frac{1}{c^2\pi^2}\|\theta_1\|_2^2\,|y|^2,$$

$$\|v\|_2\,\|\delta\|_2 \le \frac{c^2\pi^2}{4}\|v\|_2^2 + \frac{1}{c^2\pi^2}\|\delta\|_2^2$$

we get from (4.16) for all $w=(v,\varphi)\in S=\left(H^2(0,1)\cap H_0^1(0,1)\right)\times H_0^1(0,1)$, $y\in\mathbb{R}$, $\delta,\theta_1\in L^2(0,1)$:

$$\begin{aligned}\dot{\Phi}(w) &\le -\|\varphi\|_2^2 - 2\|v_x\|_2^2 + 2\frac{(c^2+1)^2\pi^2+\sigma^2}{\sigma^2c^4\pi^2}\|\theta_1\|_2^2\,|y|^2 \\ &+\pi^2\|v\|_2^2 + 2\frac{(c^2+1)^2\pi^2+\sigma^2}{\sigma^2c^4\pi^2}\|\delta\|_2^2\end{aligned} \tag{4.17}$$

Exploiting Wirtinger's inequality $\pi^2\|v\|_2^2 \le \|v_x\|_2^2$ and the fact that $\|w\|=\|(v,\varphi)\|=\sqrt{\|v_x\|_2^2+\|\varphi\|_2^2}$, we obtain from (4.17) for all $w=(v,\varphi)\in S=\left(H^2(0,1)\cap H_0^1(0,1)\right)\times H_0^1(0,1)$, $y\in\mathbb{R}$, $\delta,\theta_1\in L^2(0,1)$:

$$\dot{\Phi}(w) \le -\|w\|^2 + 2\frac{(c^2+1)^2\pi^2+\sigma^2}{\sigma^2c^4\pi^2}\|\theta_1\|_2^2\,|y|^2 + 2\frac{(c^2+1)^2\pi^2+\sigma^2}{\sigma^2c^4\pi^2}\|\delta\|_2^2 \tag{4.18}$$

It follows from (4.11), (4.14) and (4.18) that we can apply Theorem 1 with

$$K_1=\frac{1}{\sigma},\; K_2=\frac{1}{\sigma}\max\left(c^2+1+\frac{2\sigma^2}{c^2\pi^2},1+\frac{2}{c^2}\right),$$

$$X=H_0^1(0,1)\times L^2(0,1) \text{ with norm } \|w\|=\|(v,\varphi)\|=\sqrt{\|v_x\|_2^2+\|\varphi\|_2^2},$$

$$\mho=\mathcal{E}=L^2(0,1),\; S=\left(H^2(0,1)\cap H_0^1(0,1)\right)\times H_0^1(0,1) \text{ and}$$

$$p_2=1,\; R=\sqrt{2\frac{(c^2+1)^2\pi^2+\sigma^2}{\sigma^2c^4\pi^2}},\; G=2\frac{(c^2+1)^2\pi^2+\sigma^2}{\sigma^2c^4\pi^2}$$

Let $\varepsilon,\Gamma,\kappa,a,C>0$ be given (arbitrary) constants with $2\kappa>a$. The DADS feedback law is given by the equations (3.6), (3.7) where $P_i\in C^\infty(\mathbb{R})$, $i=1,2,3$ are smooth functions that satisfy inequalities (3.8), (3.9), (3.10) for all $y\in\mathbb{R}$. Applying Theorem 1 we obtain the following result for system (4.12).

**Corollary 4:** *Suppose that inequalities (4.7) hold. Let $\varepsilon,\Gamma,\kappa,a,C>0$ be constants with $2\kappa>a$. Then every solution*

$$w=(v,\varphi)\in C^0\left(\mathbb{R}_+;H_0^1(0,1)\times L^2(0,1)\right) \textit{ with } \mathbb{R}_+\ni t\to\Phi(w[t]) \textit{ being locally absolutely continuous}$$

$$\textit{and } w=(v,\varphi)\in L^1\left((0,T);\left(H^2(0,1)\cap H_0^1(0,1)\right)\times H_0^1(0,1)\right) \textit{ for each } T>0$$

$$z\in C^1(\mathbb{R}_+),\, y\in C^0(\mathbb{R}_+) \textit{ being locally absolutely continuous}$$

*of the closed-loop system (4.12) with (3.6), (3.7) that corresponds to some $d\in L^\infty(\mathbb{R}_+)$, $\theta_1\in L^\infty\left(\mathbb{R}_+;L^2(0,1)\right)$, $\theta_2\in L^\infty(\mathbb{R}_+)$, $\delta\in L^\infty\left(\mathbb{R}_+;L^2(0,1)\right)$, $b\in L^\infty\left(\mathbb{R}_+;(0,+\infty)\right)$ with* $\inf_{t\geq 0}\left(b(t)\right)>0$, *i.e., for which (4.12), (3.6), (3.7) hold for* $t\geq 0$ *a.e., satisfies estimates (3.11), (3.12), (3.13), (3.14) for all* $t\geq 0$. *Finally, if* $\lim_{t\to+\infty}\left(d(t)\right)=0$, $\lim_{t\to+\infty}\left(\delta[t]\right)=0$ *and* $\liminf_{t\to+\infty}\left(b(t)\right)\geq 1/\kappa$ *then one of the following holds:*

$$\text{either } \lim_{t\to+\infty}\left(|y(t)|\right)=\lim_{t\to+\infty}\left(\|v_x[t]\|_2\right)=\lim_{t\to+\infty}\left(\|\varphi[t]\|_2\right)=0 \text{ or}$$

$$\kappa+\exp\left(\lim_{s\to+\infty}\left(z(s)\right)\right)<\limsup_{t\to+\infty}\left(\max\left(\sqrt{2\frac{(c^2+1)^2\pi^2+\sigma^2}{\sigma^2c^4\pi^2}}\|\theta_1[t]\|_2,|\theta_2(t)|\right)\right)$$

## 5. Proof of Theorem 1

In this section we provide the proof of Theorem 1.

**Proof of Theorem 1:** In what follows we use the following inequality that is valid for all $A\geq 0, b>0$, $z\in\mathbb{R}$:

$$A\leq b\left(\kappa+\exp(z)\right)A+b\frac{\left(\kappa+\exp(z)\right)}{4a}A^2+\frac{ab}{\kappa+\exp(z)}\left(\left(\frac{1}{b}-\kappa-\exp(z)\right)^+\right)^2 \tag{5.1}$$

Inequality (5.1) follows from the facts that $A=b\frac{1}{b}A$, $\frac{1}{b}\leq\left(\kappa+\exp(z)\right)+\left(\frac{1}{b}-\kappa-\exp(z)\right)^+$ and $\left(\frac{1}{b}-\kappa-\exp(z)\right)^+A\leq\frac{\left(\kappa+\exp(z)\right)}{4a}A^2+\frac{a}{\kappa+\exp(z)}\left(\left(\frac{1}{b}-\kappa-\exp(z)\right)^+\right)^2$.

We start by proving certain inequalities for $\dot{\Phi}$ and $\dot{V}$.

Using the following inequalities that are valid for all $y,z\in\mathbb{R}$, $w\in S$, $\theta_1\in\mho$

$$R^2\|\theta_1\|_{\mho}^2\leq 2\left(\left(R\|\theta_1\|_{\mho}-\kappa-\exp(z)\right)^+\right)^2+2\left(\kappa+\exp(z)\right)^2,$$

$$2\left(\left(R\|\theta_1\|_{\mho}-\kappa-\exp(z)\right)^+\right)^2y^2\leq\frac{a}{\kappa+\exp(z)}\left(\left(R\|\theta_1\|_{\mho}-\kappa-\exp(z)\right)^+\right)^4+\frac{\kappa+\exp(z)}{a}y^4$$

we get from (3.3) for all $y,z\in\mathbb{R}$, $w\in S$, $\theta_1\in\mho$, $\delta\in\mathcal{E}$:

$$\begin{aligned}\dot{\Phi}&\leq-\|w\|^2+\frac{\kappa+\exp(z)}{a}y^4+2\left(\kappa+\exp(z)\right)^2y^2\\&+\frac{a}{\kappa+\exp(z)}\left(\left(R\|\theta_1\|_{\mho}-\kappa-\exp(z)\right)^+\right)^4+G\|\delta\|_{\mathcal{E}}^2\end{aligned} \tag{5.2}$$

Applying (5.1) with $A=\frac{\kappa+\exp(z)}{a}y^4$ and $A=2\left(\kappa+\exp(z)\right)^2 y^2$ we get from (5.2) for all $y,z\in\mathbb{R}$, $w\in S$, $\theta_1\in\mho$, $b>0$, $\delta\in\mathcal{E}$:

$$\begin{aligned}\dot{\Phi} &\leq -\|w\|^2+2b\left(\kappa+\exp(z)\right)^3 y^2\\ &+b\frac{\left(\kappa+\exp(z)\right)^5}{a}y^4\left(\kappa^{-3}+1\right)+b\frac{\left(\kappa+\exp(z)\right)^3}{4a^3}y^8\\ &+\frac{2ab}{\kappa+\exp(z)}\left(\left(\frac{1}{b}-\kappa-\exp(z)\right)^+\right)^2\\ &+\frac{a}{\kappa+\exp(z)}\left(\left(R\|\theta_1\|_{\mho}-\kappa-\exp(z)\right)^+\right)^4+G\|\delta\|_{\mathcal{E}}^2\end{aligned} \tag{5.3}$$

We also get from (3.2), (3.3) for all $y\in\mathbb{R}$, $w\in S$, $\theta_1\in\mathbb{R}^{p_1}$, $\delta\in\mathcal{E}$:

$$\dot{\Phi}\leq -\frac{1}{K_2}\Phi\left(w\right)+R^2\|\theta_1\|_{\mho}^2 y^2+G\|\delta\|_{\mathcal{E}}^2 \tag{5.4}$$

Using (3.1), (3.4), (3.5) and the fact $|\theta_2|\leq\left(|\theta_2|-\kappa-\exp(z)\right)^+ +\left(\kappa+\exp(z)\right)$, we get for all $w\in S$, $y\in\mathbb{R}$, $z\in\mathbb{R}$, $d\in\mathbb{R}$, $b>0$, $\theta_2\in\mathbb{R}^{p_2}$:

$$\begin{aligned}\dot{V} &= byu+y\left(L(w,y)\right)'\theta_2+yd\\ &\leq byu+|\theta_2||y|\left(\|w\|+\phi\left(y\right)|y|\right)+yd\\ &\leq byu+\left(|\theta_2|-\kappa-\exp(z)\right)^+|y|\|w\|+\left(\kappa+\exp(z)\right)|y|\|w\|\\ &+\left(\kappa+\exp(z)\right)\phi\left(y\right)y^2+\left(|\theta_2|-\kappa-\exp(z)\right)^+\phi\left(y\right)y^2+yd\end{aligned} \tag{5.5}$$

Exploiting the inequalities

$$yd\leq\frac{ad^2}{\kappa+\exp(z)}+\frac{\kappa+\exp(z)}{4a}y^2,$$

$$\left(\kappa+\exp(z)\right)|y|\|w\|\leq\frac{a}{4\left(\kappa+\exp(z)\right)}\|w\|^2+\frac{\left(\kappa+\exp(z)\right)^3}{a}y^2,$$

$$\begin{aligned}&\left(|\theta_2|-\kappa-\exp(z)\right)^+|y|\|w\|\\ &\leq\frac{a}{4\left(\kappa+\exp(z)\right)}\|w\|^2+\frac{\left(\kappa+\exp(z)\right)}{a}\left(\left(|\theta_2|-\kappa-\exp(z)\right)^+\right)^2y^2\end{aligned}$$

$$\frac{(\kappa+\exp(z))}{a}\left(\left(|\theta_2|-\kappa-\exp(z)\right)^+\right)^2 y^2$$
$$\leq a\frac{\left(\left(|\theta_2|-\kappa-\exp(z)\right)^+\right)^4}{\kappa+\exp(z)}+\frac{(\kappa+\exp(z))^3}{4a^3}y^4$$

$$\left(|\theta_2|-\kappa-\exp(z)\right)^+\phi(y)\,y^2\leq a\frac{\left(\left(|\theta_2|-\kappa-\exp(z)\right)^+\right)^2}{\kappa+\exp(z)}+\frac{\kappa+\exp(z)}{4a}\phi^2(y)\,y^4$$

we obtain from (5.5) for all $w\in S$, $y\in\mathbb{R}$, $z\in\mathbb{R}$, $d\in\mathbb{R}$, $b>0$, $\theta_2\in\mathbb{R}^{p_2}$:

$$\begin{aligned}\dot V\leq{}&-Cy^2+byu+\frac{(\kappa+\exp(z))^3}{4a\kappa^3}\left(\kappa+4aC+4\kappa^3+4a\kappa\phi(y)\right)y^2\\&+\frac{(\kappa+\exp(z))^3}{4a^3\kappa^2}\left(\kappa^2+a^2\phi^2(y)\right)y^4\\&+\frac{a}{\kappa+\exp(z)}\left(d^2+\frac{1}{2}\|w\|^2+\left(\left(|\theta_2|-\kappa-\exp(z)\right)^+\right)^2+\left(\left(|\theta_2|-\kappa-\exp(z)\right)^+\right)^4\right)\end{aligned}\tag{5.6}$$

Applying (5.1) with $A=\frac{(\kappa+\exp(z))^3}{4a\kappa^3}\left(\kappa+4aC+4\kappa^3+4a\kappa\phi(y)\right)y^2$ and $A=\frac{(\kappa+\exp(z))^3}{4a^3\kappa^2}\left(\kappa^2+a^2\phi^2(y)\right)y^4$ we get from (5.6) for all $w\in S$, $y\in\mathbb{R}$, $z\in\mathbb{R}$, $d\in\mathbb{R}$, $b>0$, $\theta_1\in\mho$, $\theta_2\in\mathbb{R}^{p_2}$, $\delta\in\mathcal{E}$:

$$\begin{aligned}\dot V\leq{}&-Cy^2+byu+b\frac{(\kappa+\exp(z))^7}{4a\kappa^6}\left(\kappa+4aC+4\kappa^3+4a\kappa\phi(y)\right)y^2\\&+b\frac{(\kappa+\exp(z))^7}{4a^3\kappa^5}\left(\frac{1}{16\kappa}\left(\kappa+4aC+4\kappa^3+4a\kappa\phi(y)\right)^2+\kappa^2+a^2\phi^2(y)\right)y^4\\&+b\frac{(\kappa+\exp(z))^7}{64a^7\kappa^4}\left(\kappa^2+a^2\phi^2(y)\right)^2y^8+\frac{2ab}{\kappa+\exp(z)}\left(\left(\frac{1}{b}-\kappa-\exp(z)\right)^+\right)^2\\&+\frac{a}{\kappa+\exp(z)}\left(d^2+\frac{1}{2}\|w\|^2+\left(\left(|\theta_2|-\kappa-\exp(z)\right)^+\right)^2+\left(\left(|\theta_2|-\kappa-\exp(z)\right)^+\right)^4\right)\end{aligned}\tag{5.7}$$

$$
\begin{aligned}
\dot{\Phi}+\dot{V} \leq & -\frac{2\kappa-a}{2\kappa}\|w\|^2 - Cy^2 + b\frac{(\kappa+\exp(z))^7}{4a\kappa^6}\left(\kappa+4aC+4\kappa^3+4a\kappa\phi(y)+8a\kappa^2\right)y^2 \\
& +b\frac{(\kappa+\exp(z))^7}{4a^3\kappa^5}\left(\frac{1}{16\kappa}\left(\kappa+4aC+4\kappa^3+4a\kappa\phi(y)\right)^2+\kappa^2+a^2\phi^2(y)+4a^2\left(\kappa^3+1\right)\right)y^4 \\
& +b\frac{(\kappa+\exp(z))^7}{64a^7\kappa^4}\left(\left(\kappa^2+a^2\phi^2(y)\right)^2+16a^4\right)y^8+byu+G\|\delta\|_{\mathcal{E}}^2 \\
& +\frac{a}{\kappa+\exp(z)}\left(d^2+\left(\left(|\theta_2|-\kappa-\exp(z)\right)^+\right)^2+\left(\left(|\theta_2|-\kappa-\exp(z)\right)^+\right)^4\right) \\
& +\frac{a}{\kappa+\exp(z)}\left(4b\left(\left(\frac{1}{b}-\kappa-\exp(z)\right)^+\right)^2+\left(\left(R\|\theta_1\|_{\mho}-\kappa-\exp(z)\right)^+\right)^4\right)
\end{aligned} \tag{5.8}
$$

It is clear from the fact that $2\kappa > a$ and (3.2), (3.5), (3.6), (3.8), (3.9), (3.10), (5.7), (5.8) that the following differential inequalities for all $w\in S$, $y\in\mathbb{R}$, $z\in\mathbb{R}$, $d\in\mathbb{R}$, $b>0$, $\theta_2\in\mathbb{R}^{p_2}$, $\theta_1\in\mho$, $\delta\in\mathcal{E}$:

$$
\begin{aligned}
\dot{V} \leq & -2CV(y)+\frac{a}{\kappa+\exp(z)}\left(\frac{1}{2K_1}\Phi(w)+2b\left(\left(\frac{1}{b}-\kappa-\exp(z)\right)^+\right)^2\right) \\
& +\frac{a}{\kappa+\exp(z)}\left(|d|^2+\left(\left(|\theta_2|-\kappa-\exp(z)\right)^+\right)^2+\left(\left(|\theta_2|-\kappa-\exp(z)\right)^+\right)^4\right)
\end{aligned} \tag{5.9}
$$

$$
\begin{aligned}
\dot{\Phi}+\dot{V} \leq & -\min\left(\frac{2\kappa-a}{2\kappa K_2},2C\right)\left(\Phi(w)+V(y)\right)+G\|\delta\|_{\mathcal{E}}^2 \\
& +\frac{a}{\kappa+\exp(z)}\left(|d|^2+\left(\left(|\theta_2|-\kappa-\exp(z)\right)^+\right)^2+\left(\left(|\theta_2|-\kappa-\exp(z)\right)^+\right)^4\right) \\
& +\frac{a}{\kappa+\exp(z)}\left(4b\left(\left(\frac{1}{b}-\kappa-\exp(z)\right)^+\right)^2+\left(\left(R\|\theta_1\|_{\mho}-\kappa-\exp(z)\right)^+\right)^4\right)
\end{aligned} \tag{5.10}
$$

Consider a solution $w\in C^0\left(\mathbb{R}_+;X\right)$ with $\mathbb{R}_+ \ni t\to\Phi\left(w[t]\right)$ being locally absolutely continuous and $w\in L^1\left((0,T);S\right)$ for each $T>0$, $z\in C^1\left(\mathbb{R}_+;\mathbb{R}\right)$, $y\in C^0\left(\mathbb{R}_+;\mathbb{R}\right)$ being locally absolutely continuous of the closed-loop system (3.1) with (3.6), (3.7) that corresponds to some $d\in L^\infty\left(\mathbb{R}_+\right)$, $\theta_1\in L^\infty\left(\mathbb{R}_+;\mho\right)$, $\theta_2\in L^\infty\left(\mathbb{R}_+;\mathbb{R}^{p_2}\right)$, $\delta\in L^\infty\left(\mathbb{R}_+;\mathcal{E}\right)$, $b\in L^\infty\left(\mathbb{R}_+;(0,+\infty)\right)$ with $\inf_{t\geq 0}\left(b(t)\right)>0$ .

The fact that $z(t)$ is non-decreasing (recall (3.7)) and (5.10), (3.15), (3.16), (2.16) imply the following estimate for $t\geq 0$ a.e.:

$$
\frac{d}{dt}\left(\Phi\left(w[t]\right)+V\left(y(t)\right)\right)\leq -\mu\left(\Phi\left(w[t]\right)+V\left(y(t)\right)\right)+Z \tag{5.11}
$$

Integrating the differential inequality (5.11) using the integrating factor $\exp(\mu t)$ gives the following estimate for all $t \geq 0$:

$$\Phi(w[t]) + V(y(t)) \leq \exp(-\mu t)\left(\Phi(w[0]) + V(y(0))\right) + \mu^{-1} Z \tag{5.12}$$

Estimate (3.11) is a consequence of (5.12), (3.2) and (3.5).

The fact that $z(t)$ is non-decreasing (recall (3.7)) and (5.9), (3.2), (3.5), (5.12), (3.17), (2.16) imply the following estimate for $t \geq 0$ a.e.:

$$\frac{d}{dt}\left(V(y(t))\right) \leq -2CV(y(t)) + \frac{aB}{\min(1,\kappa)\left(1+\exp(z(t))\right)} \tag{5.13}$$

Therefore, by using Lemma 4.1 in [15] we can guarantee estimate (3.12) for all $t \geq 0$.

We next show estimate (3.14). Estimate (3.12) and the fact that $z(t)$ is non-decreasing guarantees that the function $z(t)$ has a finite limit as $t \to +\infty$. This implies that the function $\exp(z(t))$ has a finite limit as $t \to +\infty$. Moreover, the facts that $\|w[\cdot]\|, y(\cdot) \in L^{\infty}(\mathbb{R}_+)$ (a consequence of (3.11)), $d(\cdot) \in L^{\infty}(\mathbb{R}_+)$, $\theta_2 \in L^{\infty}(\mathbb{R}_+;\mathbb{R}^{p_2})$, $b \in L^{\infty}(\mathbb{R}_+;(0,+\infty))$ and (3.1), (3.6), (3.4) imply that $\frac{d}{dt}(V(y(t)))$ is of class $L^{\infty}(\mathbb{R}_+)$. It follows that the function $\frac{d}{dt}(\exp(z(t))) = \Gamma\left(V(y(t))-\varepsilon\right)^{+}$ is uniformly continuous, in addition to $\int_0^{+\infty} \frac{d}{dt}(\exp(z(t)))\,dt = \lim_{t\to+\infty}\left(\exp(z(t))\right) - \exp(z(0)) < +\infty$. From Barbălat's Lemma (see [10, 18]), we have:

$$\lim_{t\to+\infty}\left(\frac{d}{dt}(\exp(z(t)))\right) = \lim_{t\to+\infty}\left(\Gamma\left(V(y(t))-\varepsilon\right)^{+}\right) = 0 \tag{5.14}$$

Therefore, by virtue of definition (3.5), estimate (3.14) holds.

We finally provide an asymptotic estimate for the unmeasured state $w$. Inequality (5.4) implies the following estimate for $t \geq 0$ a.e.:

$$\frac{d}{dt}\left(\Phi(w[t])\right) \leq -\frac{1}{K_2}\Phi(w[t]) + R^2 \|\theta_1[t]\|_{\mho}^2 y^2(t) + G\|\delta[t]\|_{\mathcal{E}}^2 \tag{5.15}$$

Integrating the differential inequality (5.15) using the integrating factor $\exp\left(\frac{1}{K_2}t\right)$ in conjunction with (3.14) gives

$$\limsup_{t\to+\infty}\left(\Phi(w[t])\right) \leq 2K_2R^2\varepsilon \limsup_{t\to+\infty}\left(\|\theta_1[t]\|_{\mho}^2\right) + K_2G\limsup_{t\to+\infty}\left(\|\delta[t]\|_{\mathcal{E}}^2\right)$$

The above estimate in conjunction with (3.2) shows that the asymptotic estimate (3.13) is valid.

Finally, combining (5.10), (3.15), (3.2) and the semigroup property we obtain:

$$\limsup_{t\to+\infty}\left(K_1\|w[t]\|^2+\frac{1}{2}|y(t)|^2\right)\le\frac{G}{\mu}\left(\limsup_{t\to+\infty}\left(\|\delta[t]\|_{\mathcal{E}}\right)\right)^2+\frac{a\left(\limsup_{t\to+\infty}\left(|d(t)|\right)\right)^2}{\mu\left(\kappa+\exp\left(\lim_{t\to+\infty}\left(z(t)\right)\right)\right)}$$
$$+\frac{a}{\mu\left(\kappa+\exp\left(\lim_{t\to+\infty}\left(z(t)\right)\right)\right)}\left(\left(\limsup_{t\to+\infty}\left(|\theta_2(t)|\right)-\kappa-\exp\left(\lim_{t\to+\infty}\left(z(t)\right)\right)\right)^+\right)^2$$
$$+\frac{a}{\mu\left(\kappa+\exp\left(\lim_{t\to+\infty}\left(z(t)\right)\right)\right)}\left(\left(\limsup_{t\to+\infty}\left(|\theta_2(t)|\right)-\kappa-\exp\left(\lim_{t\to+\infty}\left(z(t)\right)\right)\right)^+\right)^4$$
$$+\frac{4a\liminf_{t\to+\infty}\left(b(t)\right)}{\mu\left(\kappa+\exp\left(\lim_{t\to+\infty}\left(z(t)\right)\right)\right)}\left(\left(\frac{1}{\liminf_{t\to+\infty}\left(b(t)\right)}-\kappa-\exp\left(\lim_{t\to+\infty}\left(z(t)\right)\right)\right)^+\right)^2$$
$$+\frac{a}{\mu\left(\kappa+\exp\left(\lim_{t\to+\infty}\left(z(t)\right)\right)\right)}\left(\left(\limsup_{t\to+\infty}\left(R\|\theta_1[t]\|_{\mho}\right)-\kappa-\exp\left(\lim_{t\to+\infty}\left(z(t)\right)\right)\right)^+\right)^4$$

The above asymptotic estimate shows that if $\lim_{t\to+\infty}\left(d(t)\right)=0$, $\lim_{t\to+\infty}\left(\delta(t)\right)=0$, $\liminf_{t\to+\infty}\left(b(t)\right)\ge 1/\kappa$ then we have:

$$\limsup_{t\to+\infty}\left(K_1\|w[t]\|^2+\frac{1}{2}|y(t)|^2\right)$$
$$\le\frac{a}{\mu\left(\kappa+\exp\left(\lim_{t\to+\infty}\left(z(t)\right)\right)\right)}\left(\left(\limsup_{t\to+\infty}\left(\max\left(|\theta_2(t)|,R\|\theta_1[t]\|_{\mho}\right)\right)-\kappa-\exp\left(\lim_{t\to+\infty}\left(z(t)\right)\right)\right)^+\right)^2 \tag{5.16}$$
$$+\frac{2a}{\mu\left(\kappa+\exp\left(\lim_{t\to+\infty}\left(z(t)\right)\right)\right)}\left(\left(\limsup_{t\to+\infty}\left(\max\left(|\theta_2(t)|,R\|\theta_1[t]\|_{\mho}\right)\right)-\kappa-\exp\left(\lim_{t\to+\infty}\left(z(t)\right)\right)\right)^+\right)^4$$

If $\limsup_{t\to+\infty}\left(K_1\|w[t]\|^2+\frac{1}{2}|y(t)|^2\right)>0$ then $\kappa+\exp\left(\lim_{s\to+\infty}\left(z(s)\right)\right)<\limsup_{t\to+\infty}\left(\max\left(|\theta_2(t)|,R\|\theta_1[t]\|_{\mho}\right)\right)$. Otherwise, $\limsup_{t\to+\infty}\left(K_1\|w[t]\|^2+\frac{1}{2}|y(t)|^2\right)=0$ which implies that $\lim_{t\to+\infty}\left(|y(t)|\right)=\lim_{t\to+\infty}\left(\|w[t]\|\right)=0$.

The proof is complete. $\lhd$

## 6. Concluding Remarks

The present paper provided results for the partial-state regulation problem of a scalar ODE interconnected with a possibly infinite-dimensional system. In such a case the DADS control scheme allows an escape from the requirements of the small-gain theorem that is mainly used for partial-state feedback. It is clear that an extension to the case where we have more than one ODEs is needed. The case where the ODEs constitute a control system in strict feedback form is particularly interesting and is a topic for future research.